\documentclass[12pt]{article}

\usepackage[cp1251]{inputenc}
\usepackage[english]{babel}

\usepackage{amsfonts,epsfig,epstopdf}
\usepackage{amssymb,amsmath,enumerate,float}
\usepackage[section]{placeins}
\makeatother

\addto\captionsenglish{}

\renewcommand{\i}{i}
\newcommand{\ptl}{\partial}
\newcommand{\vph}{\theta}

\newcommand{\be}{\begin{equation}}
\newcommand{\ee}{\end{equation}}
\newcommand{\beq}{\begin{equation}}
\newcommand{\eeq}{\end{equation}}

\newcommand{\gH}{\mathbf{H}}
\newcommand{\gR}{\mathbf{R}}

\newcommand{\gP}{\mathbf{P}}
\newcommand{\gS}{\mathbf{S}}
\newcommand{\gC}{\mathbf{C}}

\title{Sommerfeld--type integrals for discrete diffraction problems}
\author{A.~V.~Shanin, A.~I.~Korolkov}

\begin{document}
\maketitle
\begin{abstract}
Three problems for a discrete analogue of the Helmholtz equation are studied analytically
 using the plane wave decomposition and the  Sommerfeld integral approach. 
They are: 
1) the problem with a point source on an entire plane; 
2) the problem of diffraction by a Dirichlet half-line; 
3) the problem of diffraction by a Dirichlet right angle.  
It is shown that total field can be represented as an integral of an 
algebraic function over a contour drawn on some manifold. 
The latter is a torus. As the result, the  explicit solutions are obtained  in terms of
recursive relations (for the Green's function), 
algebraic functions (for the half-line problem), or  elliptic functions
(for the right angle problem). 
 
\end{abstract}


\section*{NOTATION}
\begin{tabular}{cp{1\textwidth}}
  $\overline{\mathbb{C}}$ & Riemann sphere\\
  $K$ & wavenumber parameter \\
  $\gS$ & two-dimensional discrete lattice \\
  $\gS_2$ & 2-sheet branched discrete lattice\\
  $\gS_3$ & 3-sheet branched discrete lattice\\
  $u(m,n)$ & wave field on a lattice $\gS$\\
  $U(\xi_1, \xi_2)$ & Fourier transform of $u$\\
  $(r , \phi)$ & polar coordinates on the plane $(m,n)$ \\
  $\phi_{\rm in}$ & angle of propagation of the incident wave \\
  $x = e^{i \xi_1}$, $y = e^{i \xi_2}$ & ``algebraic wavenumbers'' \\ 
  $\tilde u(m,n)$ & total field on a discrete branched surface \\
  $\gC$ & Circle $[0,2\pi]$ with $0$ and $2\pi$ glued to each other \\ 
  $ D(\xi_1, \xi_2)$ & dispersion function, defined by (\ref{eq0105}) \\
  $\hat D(x, y)$ & dispersion function, defined by (\ref{eq0105a}) \\
  $\Xi(x)$ & root of dispersion equation (\ref{eq0111}), defined by (\ref{eq0111b})\\
  $\gR$ & Riemann surface of $\Xi(x)$ \\
  $\gR_2$ & 2-sheet covering of $\gR$\\
  $\gR_3$ & 3-sheet covering of $\gR$\\  
  $\gH$ & dispersion surface, the set of all points $(x,y),x,y \in\overline{\mathbb{C}}$ such that (\ref{eq0111}) is valid \\
  $\gH_2$ & 2-sheet covering of $\gH$\\
  $\gH_3$ & 3-sheet covering of $\gH$\\ 
  $B_1 \dots B_4$ & branch points of $\gR$ \\
  $J_1 \dots J_4$ & zero / infinity points on $\gH$\\ 
  $\sigma_1 \dots \sigma_4$ &  integration contours for the representations 
  	(\ref{eq0107b}), (\ref{eq0107g}), (\ref{eq0107h}), (\ref{eq0107j}) \\
  $\Gamma_j$ & contours for the Sommerfeld integral \\	 						
  $(\alpha, \beta) \in \gC^2$ & coordinates on the torus\\  
  $\zeta$ & mapping between $\gH$ and $\gR$\\
  $\Psi$ & analytic 1-form  on $\gH$ defined by (\ref{eq0301}) \\
  $f_0(x) , \dots ,f_3(x)$ & basis algebraic functions 
  						on $\gR_2$  \\
  $t(p)$ & elliptic integral (\ref{elfunF}) \\
  $E(t, \omega_1, \omega_2)$ & elliptic function (\ref{etafunct})  						 
\end{tabular}

\section{Introduction}
In the beginning of 20th century, Sommerfeld found a closed integral solution for the 
problem of diffraction by a half-plane \cite{Sommerfeld1954}
by combining the plane wave decomposition integral with the 
reflection method. Later on, this plane wave decomposition integral with some particular contour of integration 
has been named after Sommerfeld. 
The Sommerfeld integral approach has been then applied to a number of problems such as problem of diffraction by  a wedge \cite{Sommerfeld1901,MacDonald1902,Babich2008} and some others \cite{Luneburg1997,Hannay2003}.

In this paper, we build an analogy of the Sommerfeld integral for discrete problems.  The following problems on a 2D square lattice are considered in the 
paper:
\begin{itemize}
\item radiation of a point source in an entire plane 
(i.~e.,\ finding of a Green's function of a plane),
\item diffraction by a Dirichlet half-line,
\item diffraction by a Dirichlet right angle.
\end{itemize} 

The discrete Green's function problem has been studied 
in two different physical contexts. They are the discrete potential theory \cite{Duffin1953, Thomassen1990, Zemanian1988, Atkinson1999},  and the problem of the random walk \cite{Pearson1905, McCrea1940, Spitzer1964}. Both problems can be reduced to the calculation of the Green's function for a discrete Laplace equation \cite{Ito1960}. The Green's 
function in these works is represented in terms of a double Fourier integral. Depending on the position of the observation point, four different single integral representations can be introduced with the help of residue integration. In the current work we study the representation as an integral on the complex manifold of dimension~1. The latter is a torus.  

The problem of diffraction by a half-line is well known in the context of fracture mechanics \cite{Slepyan1982, Slepyan2002}, where the Dirichlet 
half-line models a rigid constraint in a square lattice. The problem has been 
solved by several authors \cite{Eatwell1982, Slepyan1982, Sharma2015b} 
with the help of the Wiener-Hopf technique. 
The solution has been expressed in terms of elliptic integrals. Here we introduce an analogue of the Sommerfeld integral for this problem and obtain an expression for the solution in terms of algebraic functions. 

We are not aware of any analytical results on right angle diffraction problem
for a lattice. Here we obtain a solution of this problem in terms of elliptic functions 
by using the Sommerfeld integral approach.    


\section{Discrete Green's function on a plane}

\subsection{Problem formulation}
\label{Sec:2.1}

Let there exist a two-dimensional lattice, whose nodes are indexed by  
$m,n \in \mathbb{Z}$. This lattice is referred to as $\gS$. 
Let a function $u(m,n)$ defined on $\gS$ obey the equation
\begin{equation}
u(m+1, n) + u(m-1,n) + u(m,n-1) + u(m,n+1) + (K^2 - 4) u(m,n) = 
\delta_{m,0} \delta_{n,0},
\label{eq0101}
\end{equation}
where $\delta_{m,n}$ is the Kronecker's delta.
Indeed the expression 
\[
u(m+1, n) + u(m-1,n) + u(m,n-1) + u(m,n+1) - 4\, u(m,n) 
\] 
is the discrete analogue of the continuous 2D Laplace operator $\Delta$.
Our aim is to compute $u(m,n)$.

The wavenumber parameter $K$ is close to positive real, but has a small 
positive imaginary part mimicking an attenuation on the lattice. 
The radiation condition (in the form of the limiting absorption principle) 
states that $u(m,n)$ should decay exponentially as $\sqrt{m^2 + n^2} \to \infty$.


\subsection{Double integral representation. Dispersion equation}
\label{Sec:2.2}
Apply a double Fourier transform to (\ref{eq0101}). The transform is as follows:  
\begin{equation}
u(m,n) 
\, 
\longrightarrow
\,
U(\xi_1, \xi_2) = \sum_{m,n = -\infty}^{\infty}
u(m,n) \exp \{ - i (m \xi_1 + n \xi_2)  \} .
\label{eq0102}
\end{equation} 
The result is 
\begin{equation}
D(\xi_1, \xi_2) U(\xi_1 , \xi_2) = 1,
\label{eq0102a}
\end{equation}
where 
\begin{equation}
D(\xi_1 , \xi_2) \equiv
2 \cos \xi_1 + 2 \cos \xi_2 - 4 + K^2.
\label{eq0105}
\end{equation}
The inverse Fourier transform is given by 
\begin{equation}
U(\xi_1 , \xi_2)
\,
\longrightarrow
\,
u(m, n) = \frac{1}{4 \pi^2} 
\int \! \! \! \! \int_{-\pi} ^{\pi}
U(\xi_1,\xi_2) \exp \{  i (m \xi_1 + n \xi_2)  \}\, d\xi_1 \, d\xi_2 ,  
\label{eq0103}
\end{equation} 
and thus the following representation for $u(m,n)$ holds: 
\begin{equation}
u(m,n) = \frac{1}{4\pi^2} 
\int \! \! \! \! \int_{-\pi} ^{\pi}
\frac{
\exp \{  i (m \xi_1 + n \xi_2)  \}
}{D(\xi_1 , \xi_2)}
d\xi_1 \, d\xi_2 . 
\label{eq0104}
\end{equation}

Introduce the variables 
\begin{equation}
x = e^{i \xi_1} , 
\qquad 
y = e^{i \xi_2}, 
\label{eq0110}
\end{equation}
and the function 
\begin{equation}
\hat D(x , y) \equiv
x + x^{-1} + y + y^{-1} - 4 + K^2.
\label{eq0105a}
\end{equation}
(Indeed, $D(\xi_1, \xi_2) = \hat D(e^{i \xi_1}, e^{i\xi_2})$.)
The integral (\ref{eq0104}) can be rewritten as 
\begin{equation}
u(m,n)
= - \frac{1}{4 \pi^2}
\int_\sigma \int_\sigma 
\frac{x^m y^n}{\hat D(x , y)} \frac{dx}{x} \frac{dy}{y},
\label{eq0104a}
\end{equation}
where 
contour
$\sigma$ is the unit circle in the $x$-plane ($|x|=1$) passed 
in the positive direction (anti-clockwise).
Expression (\ref{eq0104a}) is the double integral 
representation for the Green's function $u$. 

The combination 
\[
x^m y^n = e^{i (m \xi_1 + n \xi_2)}
\]
plays the role of a plane wave on a lattice. 
If $x$ and $y$ obey the {\em dispersion equation}
\begin{equation}
\hat D (x, y) = 0,
\label{eq0111}
\end{equation}
then such a wave can travel along a lattice, being supported by 
a homogeneous equation 
\begin{equation}
u(m+1, n) + u(m-1,n) + u(m,n-1) + u(m,n+1) + (K^2 - 4) \, u(m,n) = 0.
\label{eq0101h}
\end{equation}

Among general plane waves corresponding to any solution of (\ref{eq0111}) 
we would like to select the subset of {\em real waves}. 
If $K$ is real, real waves are just waves with $|x| = |y| = 1$. 
They are  plane waves in the usual understanding (in contrast with waves 
attenuating in some direction). Such waves can be characterized by a
(real) wavenumber $(\xi_1, \xi_2)$, or by the propagation angle 
$\vph$ such that 
\begin{equation}
\xi_1 = \xi \cos \vph , 
\qquad 
\xi_2 = \xi \sin \vph , 
\label{eq4101}
\end{equation}
and $\xi = \xi(\vph)$ is real positive. 

Angle $\vph$ takes values on the circle $[0, 2\pi]$ with $0$ and $2\pi$
glued to each other. Below we refer to this circle as $\gC$.
 
It follows from (\ref{eq4101}) that 
\begin{equation}
{\rm Im} [\xi_1 / \xi_2] = 0.
\label{eq4102}
\end{equation} 
So, if $K$ is real then real waves correspond to 
solutions of dispersion equation $D(\xi_1 , \xi_2) =0$ 
obeying (\ref{eq4102}). 
Note that not all such pairs $(\xi_1, \xi_2)$
correspond to the real waves. There 
are two branches of such pairs, both organized as $\gC$, and 
only one branch corresponds to real $(\xi_1, \xi_2)$,
i.~e.\ to the real waves. 

Mainly for clarity and convenience, we are going to introduce ``real waves'' for the 
case of complex~$K$. In this case, there are no non-decaying  plane waves, so there exists 
an ambiguity in the choice of the real waves. We solve this ambiguity by selecting 
the solutions of  the dispersion equation 
with 
\begin{equation}
{\rm Im} [\sin (\xi_1)  / \sin (\xi_2)] = 0,
\label{eq4102a}
\end{equation} 
or, the same, the solutions of (\ref{eq0111}) with  
\begin{equation}
{\rm Im} \left[
\frac{x - x^{-1}}{y - y^{-1}}
\right] = 0.
\label{rewaves}
\end{equation}
One can show that there is a branch of such pairs $(x, y)$ having the shape of a loop 
on the Riemann surface of $y(x)$ that tends to usual real waves as 
${\rm Im}[K]\to 0$. Topologically, the real waves remain organized as $\gC$.

The choice of the relation (\ref{rewaves}) for the real waves is explained below.   
Namely, the saddle points of the integral representations of the field found from equation 
(\ref{eq:sadpoint}) belong to the set described by (\ref{rewaves}).

\subsection{Single integral representations}
\label{Sec:2.3}
The integral (\ref{eq0104a}) can be taken with respect to one of the variables by the 
residue integration. As the result, one can obtain a single integral representation.
There are four cases, possibly intersecting:
\[
n \ge 0, \qquad 
n \le 0, \qquad 
m \ge 0, \qquad 
m \le 0. 
\]
Each of these cases results in its own single integral representation formula. 

\begin{figure}[ht]
\centerline{\epsfig{file=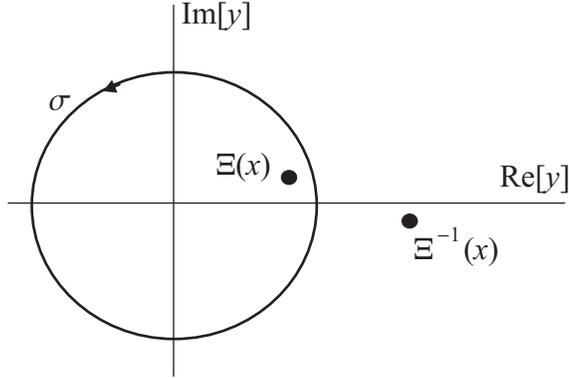}}
\caption{Complex plane $y$ for fixed $x$}
\label{fig03b}
\end{figure}

\vskip 6pt
\noindent
{\bf Case $n \ge 0$}:

Consider the integral (\ref{eq0104a}). 
Fix $x \in \sigma$ and study the integral with respect to~$y$. The $y$-plane is shown in 
Fig.~\ref{fig03b}. One can see that there are four possible singular points in this plane. 

Two of them are the roots of 
equation (\ref{eq0111}) considered with respect to $y$. 
These roots are 
\begin{equation}
y(x) = \Xi (x), \qquad y(x) = \Xi^{-1} (x),
\label{eq0111a}
\end{equation}
where 
\begin{equation}
\Xi (x) = -\frac{K^2 - 4 + x + x^{-1}}{2}
+
\frac{\sqrt{(K^2 - 4 + x + x^{-1})^2 - 4}}{2} .
\label{eq0111b}
\end{equation}

The value of the square root is chosen in such a way that 
$|\Xi (x)| < 1$. Note that $|\Xi (x)|$ cannot be equal to~1
if $|x| = 1$ since $K$  is not real. 
The points (\ref{eq0111a}) are simple poles of the 
integrand of (\ref{eq0104a}).

Note also that 
\begin{equation}
\Xi^{-1} (x) = -\frac{K^2 - 4 + x + x^{-1}}{2}
-
\frac{\sqrt{(K^2 - 4 + x + x^{-1})^2 - 4}}{2}, 
\label{eq0111p}
\end{equation}
and, indeed, $y = \Xi(x)$ and $y = \Xi^{-1}(x)$ are two roots  
of the quadratic 
equation (\ref{eq0111}).

Beside (\ref{eq0111a}), there maybe singularities 
of the integrand
at two other points: $y = 0$ and $y= \infty$, 
(the latter is a certain point of the Riemann sphere 
$\overline{\mathbb{C}}$).
The presence of singularities at these points depends on the value of~$n$.
If $n \ge 0$ (as is in the case under consideration) 
then the integrand is regular at $y = 0$ and may have a pole at $y = \infty$.  

Thus, $y = \Xi(x)$ is the only singularity of the integrand inside 
the contour $\sigma$.
Apply the residue theorem. The result is 
\begin{equation}
u(m,n) = \frac{1}{2\pi i} 
\int_\sigma
\frac{
x^m \, \Xi^n(x)
}{
\Xi(x) \,
\ptl_{y} \hat D( x , \Xi(x) ) 
}
\frac{d x}{x} .
\label{eq0107}
\end{equation} 
As one can find by direct computation, 
\begin{equation}
y \, \ptl_y \hat D (x, y) = y - y^{-1}. 
\label{eq0107a}
\end{equation}
Thus, (\ref{eq0107}) can be rewritten as 
\begin{equation}
u(m,n) = \frac{1}{2\pi i} 
\int_\sigma
x^m \, y^n
\frac{
d x
}{
x (y - y^{-1})
},
\qquad 
y = \Xi(x).
\label{eq0107b}
\end{equation} 
This is a single integral representation of the field. 

\vskip 6pt
\noindent
{\bf Case $n \le 0$}: 

In this case $y = 0$ is a singular point of the integrand of (\ref{eq0107}), 
but $y = \infty$ is a regular point (more rigorously, 
$y = \infty$ is a regular point of the differential form
that is integrated 
in (\ref{eq0107})). 
This means that the integrand has no branching at $y = \infty$, 
and it decays not slower than $\sim y^{-2}$. For such an integrand, 
one can apply the residue theorem to the 
exterior of $\sigma$. The result is 
\begin{equation}
u(m,n) = -\frac{1}{2\pi i} 
\int_\sigma
x^m \, y^n
\frac{
d x
}{
x (y - y^{-1})
},
\qquad 
y = \Xi^{-1}(x).
\label{eq0107g}
\end{equation} 

\vskip 6pt
\noindent
{\bf Case $m \ge 0$}:

The representation of the field is 
\begin{equation}
u(m,n) = \frac{1}{2\pi i} 
\int_\sigma
x^m \, y^n
\frac{
d y
}{
y (x - x^{-1})
},
\qquad 
x = \Xi(y) .
\label{eq0107h}
\end{equation} 

\vskip 6pt
\noindent
{\bf Case $m \le 0$}: 

The representation of the field is 
\begin{equation}
u(m,n) = -\frac{1}{2\pi i} 
\int_\sigma
x^m \, y^n
\frac{
d y
}{
y (x - x^{-1})
},
\qquad 
x = \Xi^{-1}(y).
\label{eq0107j}
\end{equation} 


\vskip 6pt

Thus, we obtained four single integral representations: (\ref{eq0107b}), 
(\ref{eq0107g}), (\ref{eq0107h}), (\ref{eq0107j}).



\subsection{Field representation by integration on a manifold}
\label{Sec:2.4}
Let us analyze the integral (\ref{eq0107b}).  
Consider $x$ and 
$y$ as  complex variables taking values on the Riemann sphere  $\mathbb{\bar C}$.
We remind that $\mathbb{\bar  C}$ is a compactified complex plane, i.~e.\  a plane to which 
the infinite point is added. The usage of the Riemann sphere is convenient 
when it is necessary to study functions having algebraic growth, or just algebraic functions,
which is the case in (\ref{eq0107b}).

Each point $(x, y)$ thus belongs to $\overline{\mathbb{C}} \times 
\overline{\mathbb{C}}$. 
Let us describe
the set of points $(x, y)$ such that  equation (\ref{eq0111}) is valid. 
It is easy to prove that this set is an analytic manifold of complex dimension~1 
(so it has real dimension~2). We comment this below. 
This manifold will be referred to as $\gH$.
Since (\ref{eq0111}) is the dispersion relation for the lattice, one can call 
$\gH$ the {\em dispersion surface}.

The manifold $\gH$ can be easily built using the function $\Xi(x)$  defined by (\ref{eq0111b}).
This function has been defined only for $x \in \sigma$. Continue this 
function analytically onto the whole 
$\overline{\mathbb{C}}$. 
This function becomes double-valued, with some branch points. 
Indeed, $\gH$ is the set of all points $(x , \Xi(x))$,
$x \in \overline{\mathbb{C}}$. 

We make here an obvious observation that the Riemann surface of $\Xi(x)$, 
which will be referred to as $\gR$,  
is the projection of $\gH$ onto~$x$. Thus, topologically $\gH$
coincides with the Riemann surface $\gR$. 
We will denote the mapping between $\gH$ and $\gR$ by $\zeta$.
More often, we will use the inverse mapping 
\[
x \stackrel{\zeta^{-1}}{\longrightarrow} (x , \Xi(x))
,
\qquad
x \in \gR.
\]

Let us study the Riemann surface $\gR$.
Function $\Xi(x)$ has four  branch points. They are the points where the argument of the 
square root in (\ref{eq0111b}) is equal to zero. 
By solving the equation 
\[
(K^2 - 4 + x + x^{-1})^2 - 4 =0
\] 
we find that the branch points are $x = \eta_{j,k}$, $j, k = 1,2$, 
where 
\begin{equation}
\eta_{1,1} = -\frac{d}{2} + \frac{\sqrt{d^2 -4}}{2},
\qquad d = K^2 - 2 ,  
\label{eq0111c1}
\end{equation} 
\begin{equation}
\eta_{1,2} = -\frac{d}{2} + \frac{\sqrt{d^2 -4}}{2},
\qquad d = K^2 - 6 ,  
\label{eq0111c2}
\end{equation} 
\begin{equation}
\eta_{2,1} = -\frac{d}{2} - \frac{\sqrt{d^2 -4}}{2},
\qquad d = K^2 - 2 ,  
\label{eq0111c3}
\end{equation} 
\begin{equation}
\eta_{2,2} = -\frac{d}{2} - \frac{\sqrt{d^2 -4}}{2},
\qquad d = K^2 - 6 .  
\label{eq0111c4}
\end{equation} 

Let us list some important properties of the branch points that can be checked directly 
or derived from elementary properties of quadratic equations. These properties are: 
\begin{itemize}
\item 
The branch points are the points at which $\Xi(x) = \pm 1$. 

\item 
For $y = \Xi(x)$
\begin{equation}
\Upsilon(x) \equiv x (y - y^{-1}) = \sqrt{(x - \eta_{1,1})(x - \eta_{1,2})(x - \eta_{2,1})(x - \eta_{2,2})}. 
\label{eq0111w}
\end{equation}
Note that the left-hand side of (\ref{eq0111w}) is the denominator of the integrand of 
(\ref{eq0107b}).  
\item 
$
\eta_{1,1} \,  \eta_{2,1} = 1, 
\qquad  
\eta_{1,2} \, \eta_{2,2} = 1. 
$

\item
Exactly two of the branch points
$(\eta_{1,1}, \eta_{1,2}, \eta_{2,1} , \eta_{2,2})$
are located inside the circle $|x|<1$.
By the choice of the square root branches, we can make these points 
be called $\eta_{2,1}$ and $\eta_{2,2}$. 

\end{itemize}

\begin{figure}[ht]
\centerline{\epsfig{file=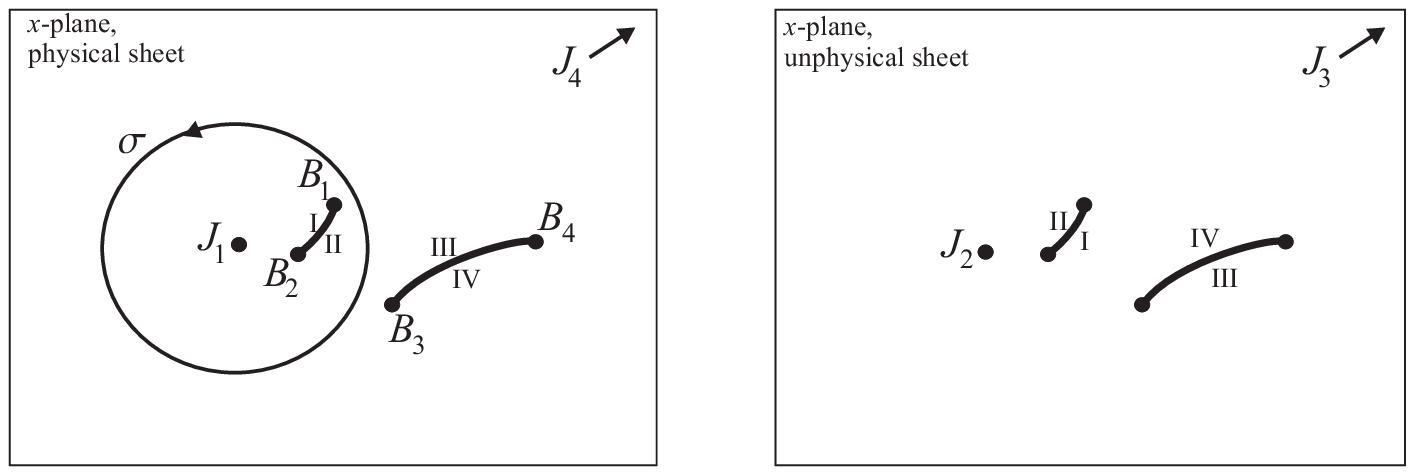}}
\caption{Scheme of $\gR$}
\label{fig03}
\end{figure}

The scheme of $\gR$ is shown in Fig.~\ref{fig03}.
Two sheets $\gR$ are Riemann spheres.
They are shown projected onto a plane, such that the infinity point 
of $\overline{\mathbb{C}}$ becomes infinitely remote.

The branch points are connected by cuts shown by bold curves. For definiteness, the branch cuts are conducted along the lines at which $|\Xi(x)|=1$.
The shores of the cuts labeled by equal Roman number should be connected with each other. 

One of the sheets drawn in Fig.~\ref{fig03} is called {\em physical\/}, and
the other is {\em unphysical\/}.
The physical sheet is the one on which $|\Xi (x)| < 1$ for $|x| = 1$. Respectively, 
on the unphysical sheet 
$|\Xi (x)| > 1$ for $|x| = 1$.
The integration in (\ref{eq0107b}) is taken along contour $\sigma$
drawn on the physical sheet. 

It is convenient to label the points of $\gR$ by 
coordinates of corresponding points of $\gH$, i.~e.\ by the pairs 
$(x, \Xi(x))$. 

There are several important points on this Riemann surface. 
First, they are the branch points 
\[ 
B_1\,:\,(\eta_{2,1},1),
\qquad  
B_2\,:\,(\eta_{2,2},-1),
\qquad  
B_3\,:\,(\eta_{1,1},1),
\qquad  
B_4\,:\,(\eta_{1,2},-1).
\]
Second, they are zero/infinity points, i.~e.\ the points at which either 
$x$ or $y = \Xi(z)$ is zero or infinity. They are the points 
\[
J_1 \, : \, (0,0), 
\qquad 
J_2 \, : \, (0,\infty), 
\qquad 
J_3 \, : \, (\infty, \infty), 
\qquad 
J_4 \, : \, (\infty,0). 
\]

\begin{figure}[ht]
\centerline{\epsfig{file=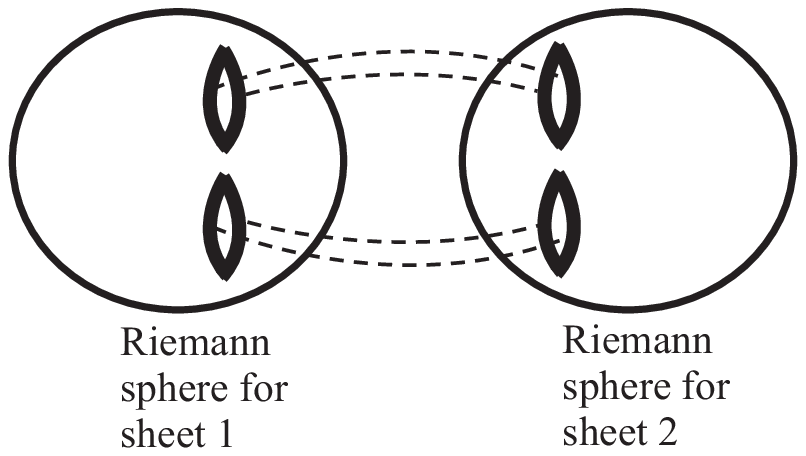, width = 6cm}}
\caption{Why $\gR$ is a torus}
\label{fig03d}
\end{figure}

Topologically, $\gH$ is a torus (i.~e.\ it has genus equal to 1). This can be easily understood, since $\gH$ is obtained by taking two spheres, making two cuts, and connecting their shores. The scheme of making a torus out of two Riemann spheres is shown in 
Fig.~\ref{fig03d}.  
 
The statement that $\gH$ is an analytic manifold means that 
in each (small enough) neighborhood of any point of $\gH$ one can introduce a 
complex local variable, such that all transformation mappings between the 
neighboring local variables are biholomorphic. It is clear that such local 
variables can be: 
\begin{itemize}
\item
$x$ for all points except the branch points 
$B_1, \dots , B_4$, and 
and two 
infinities $J_3$ and $J_4$; 

\item
$y$ for the branch points $B_1, \dots , B_4$; 

\item
$\tau = 1/x$ for the infinities $J_3$ and $J_4$.  
\end{itemize}   
   
\begin{figure}[ht]
\centerline{\epsfig{file=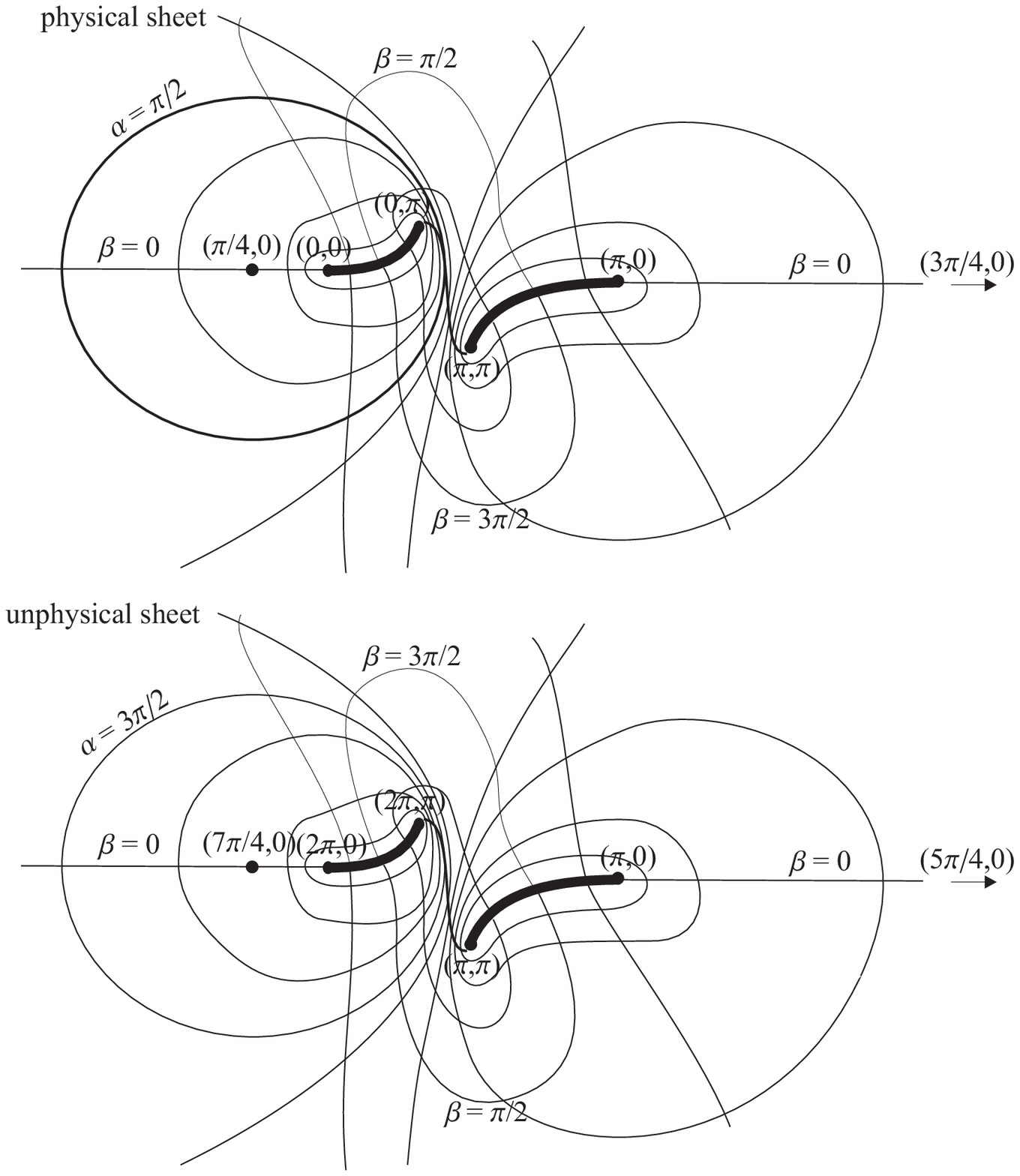}}
\caption{Coordinates $(\alpha, \beta)$ on $\gR$}
\label{fig11a}
\end{figure}

To gain some clarity, we introduce coordinates $(\alpha, \beta)$
on $\gH$ showing that this is a torus. 
Both coordinates are real and take values in~$\gC$. 
The coordinate lines on $\gH$ (projected on $\gR$) are shown in Fig.~\ref{fig11a}.
The explicit formulae for the coordinates are not important 
(for the topological purposes the coordinates 
can be drawn just ``by hands''), but we keep the following properties valid: 
\begin{itemize}
\item
Points $B_1$, $B_2$, $B_3$, $B_4$ have coordinates 
$(0, \pi)$, $(0,0)$, $(\pi, \pi)$, $(\pi, 0)$, respectively.

\item
Points $J_1$, $J_2$, $J_3$, $J_4$ have coordinates 
$(\pi/4, 0)$, $(7\pi/4,0)$, $(5\pi /4, 0)$, $(3\pi / 4, 0)$, respectively.

\item
Contour $\sigma$ corresponds to the line $\alpha = \pi/2$ passed in the 
negative direction.

\item 
The cuts (bold lines), taken for $|y| = 1$, correspond to $\alpha = 0$
and $\alpha = \pi$.

\item 
There is an important set of points on $\gH$ where relation 
(\ref{rewaves}) is fulfilled.
A study of the explicit expressions for this set shows that it consists 
of two loops. 
One of the loops passes through the points 
$B_1$ and $B_3$. This is the {\em real waves\/} line discussed above. 
We force the coordinate line 
$\beta = \pi$ coincide with this line. 
The other loop bears all infinity 
points and the branch points $B_2$ and $B_4$. We make the coordinate line 
$\beta = 0$ coincide with this loop.   

\item 
On the line $\beta = \pi$ (the real waves line) we force the coordinate $\alpha$ to have values 
\begin{equation}
\alpha = \arctan \left(
\frac{y-y^{-1}}{x - x^{-1}} 
\right) = \arctan \left(
\frac{\sin \xi_2}{\sin \xi_1} 
\right).
\label{rewaves1}
\end{equation}
Our $\arctan$ function takes values in $\gC$ (not in $[-\pi/2, \pi/2]$).  
We assume that $\alpha = 0$ at $B_2$, and then use (\ref{rewaves}) taking the 
values on $\gC$ by continuity.



\end{itemize}

Note that the variables $(\alpha, \beta)$ are real, and they are used for display purposes only. They are not related to the complex analytic 
structure on the torus (for which one needs a single complex variable). 
Some ``proper'' coordinates will be introduced on $\gH$ by using the 
elliptic variable $t$ in Section~\ref{wedge4}. 

The introduction of the manifold $\gH$ (dispersion surface) is needed mainly to 
obtain the field representation invariant to the change of variables. Beside 
$\gH$ itself, we 
need a formalism of analytic differential forms on $\gH$ and integrating of them. 
An analytic 1-form can be defined in the manifold $\gH$ (see e.g.~\cite{Gurvitz1968,Shabat1992}) by introducing 
a formal expression $h(z) dz$, where $z$ is a local complex variable
for some neighborhood, and $h(z)$ is an analytic function in this 
neighborhood. In intersecting neighborhoods the representations 
of a form
can be different, say $h_1 (z_1) dz_1$ and $h_2 (z_2) dz_2$, but they should match
in an obvious way:
\[
h_2 = h_1 \frac{d z_1}{d z_2}. 
\] 

The 1-form can be analytic/meromorphic if the functions $h_j$ are analytic/meromorphic. 
In the same sense the form can have a zero or a pole of some order.

Analyticity of a 1-form is an important property since if a form is analytic in 
some domain of an analytic manifold, and this form is integrated along some contour, 
this contour can be deformed within this domain, i.~e.\ a usual Cauchy's theorem 
for a complex plane can be generalized onto an analytic manifold. Integration over the 
poles also keeps the same. 

Let us prove that the form 
\begin{equation}
\Psi = \frac{dx}{x (y - y^{-1})},
\label{eq0301}
\end{equation}
which is a part of (\ref{eq0107g}),
is analytic everywhere on $\gH$. The statement is trivial everywhere 
except the infinities and the branch points. Consider the infinities. 
At the points $J_1$ and $J_2$ it is easy to show that $(y-y^{-1}) \sim x^{-1}$
as $x \to 0$, thus the 
denominator is non-zero. At the points $J_3$ and $J_4$ one can show that
$(y-y^{-1}) \sim x$ as $x \to \infty$, thus 
$ \Psi \sim x^{-2} dx $. A change to the variable $\tau = 1/x$ shows that the form is regular. 

Finally, consider the branch points $B_1 , \dots , B_4$. As it has been 
mentioned, one can take $y$ as a local variable at these points. 
An important observation is that due to the theorem about an implicit function,  
\begin{equation}
\frac{d y}{d x} = - \frac{\ptl_x \hat D}{\ptl_y \hat D} 
\label{eq0115}
\end{equation}
everywhere on $\gH$.
Thus, 
\begin{equation}
\frac{dx}{x(y - y^{-1})} = - \frac{dy}{y(x-x^{-1})}. 
\label{eq0116}
\end{equation}
The denominator of the right-hand side 
of (\ref{eq0116}) is non-zero at the branch points, so the form is regular. 

The representation (\ref{eq0107b}) can be rewritten as a contour integral 
of the form 
\begin{equation}
\psi_{m,n} = \frac{i}{2\pi} x^m y^n \Psi 
\label{eq0302}
\end{equation}
along some contour drawn directly on~$\gH$.
The contour is, indeed, the preimage of $\sigma$ shown in Fig.~\ref{fig03},
i.~e.\ $\zeta^{-1} (\sigma)$. 

\begin{figure}[ht]
\centerline{\epsfig{file=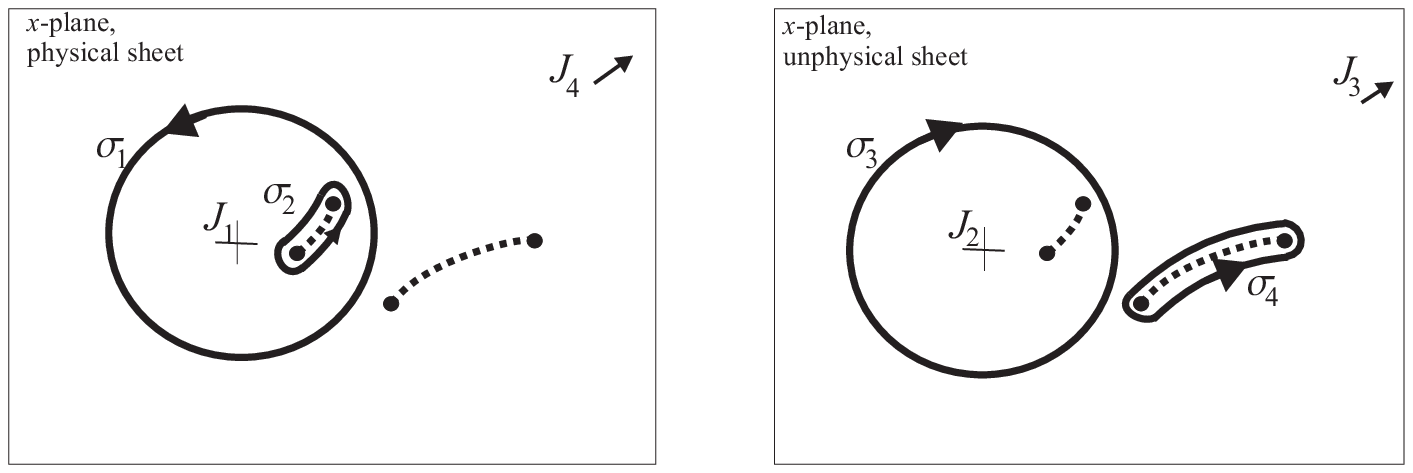}}
\caption{Position of contours and infinities on the Riemann surface}
\label{fig07}
\end{figure}

It can be easily shown that 
three other representations, (\ref{eq0107g}),
(\ref{eq0107h}), (\ref{eq0107j}), can be written as contour integrals 
of {\em the same differential form\/} $\psi_{m,n}$ on $\gH$, but taken along some other contours. 
Namely, for the integrals 
(\ref{eq0107g}),
(\ref{eq0107h}), 
(\ref{eq0107j}), 
these contours projected onto $\gR$
are shown in Fig.~\ref{fig07}. They are denoted by $\sigma_3$, $\sigma_2$, 
$\sigma_4$, respectively. 
The contour  $\sigma$ is denoted by $\sigma_1$ for uniformity.
The contours on $\gH$ are $\zeta^{-1} (\sigma_j)$ for $j = 1,\dots ,4$. 

The contours $\sigma_1$, $\sigma_2$, $\sigma_3$, $\sigma_4$ 
in the coordinates $(\alpha, \beta)$ correspond to coordinate 
lines $\alpha = \pi/2, \, 0, \, 3 \pi/2, \, \pi$, 
respectively. The direction of bypass of each contour 
is negative with respect to the variable~$\beta$.

The representations (\ref{eq0107b})--(\ref{eq0107j})
can be written in the common form 
\begin{equation}
u(m,n) = \int_{\zeta^{-1}(\sigma_j)} w_{m,n} \Psi,
\label{eq0303}
\end{equation}
where 
\begin{equation}
w_{m,n} = w_{m,n} (x,y) = x^m y^n 
\label{eq0304}
\end{equation}
is the ``plane wave''. Note that the integration is held over a contour on 
$\gH$,
so $(x,y) \in \gH$, and any such $w_{m,n}$ obey the homogeneous 
stencil equation (\ref{eq0101h}).  
Representation (\ref{eq0303}) can be considered as a 
generalized 
plane wave decomposition. 

Below we don't distinguish between $\gH$ and~$\gR$ and write $\sigma_j$ instead of $\zeta^{-1} (\sigma_j)$ if it does not lead to an ambiguity.

The form $\psi_{m,n}$ is  analytic everywhere on $\gH$
only for $m= n= 0$. 
Depending on $m$ and $n$, this form can have poles at the infinity points. 
For each of the points $J_1 , \dots , J_4$
one can derive a condition providing regularity of that point. 
These conditions of regularity for the infinity points are as follows: 
\begin{eqnarray*}
J_1 &:& m+n \ge 0,
\\
J_2 &:& m-n \ge 0,
\\
J_3 &:&-m-n \ge 0,
\\
J_4 &:&-m+n \ge 0.
\end{eqnarray*}

%
%

\begin{figure}[ht]
\centerline{\epsfig{file=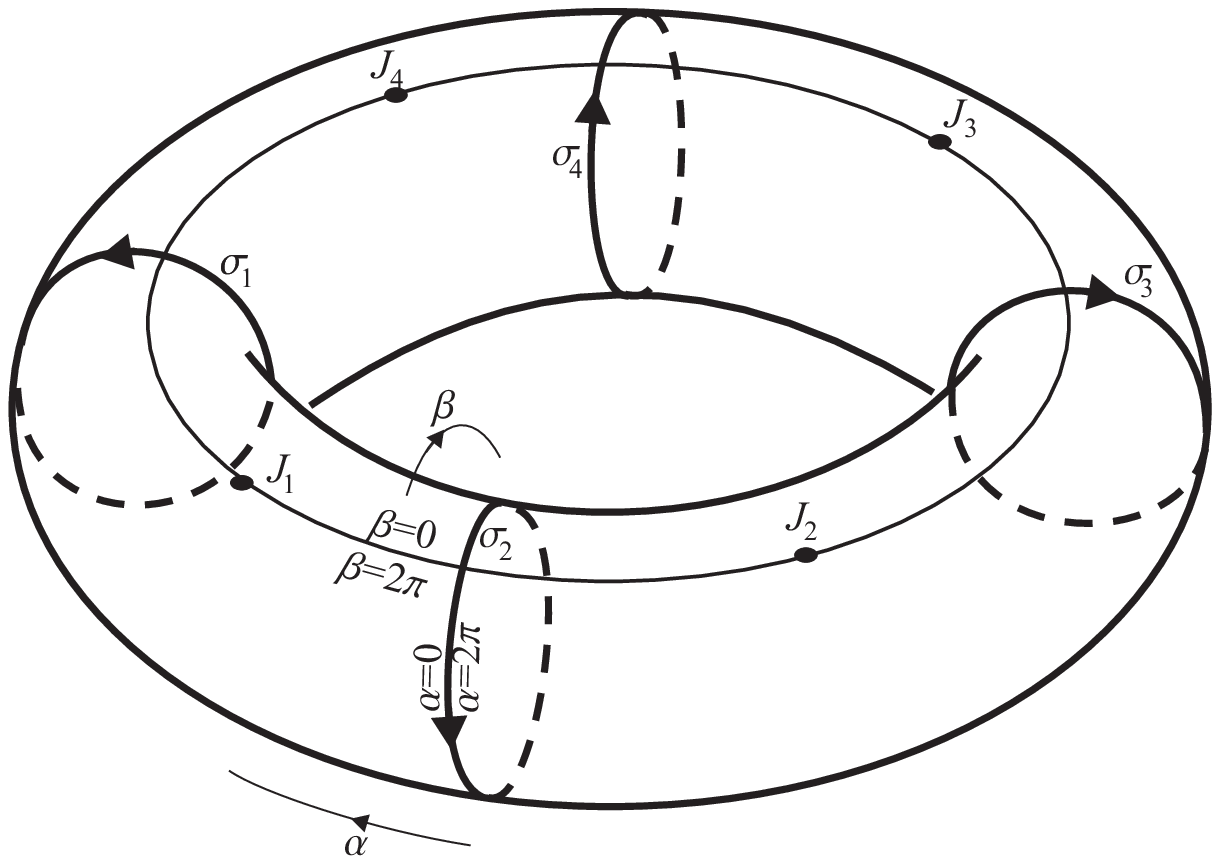,width = 9cm}}
\caption{Position of contours and infinities on $\gH$ in the coordinates $(\alpha, \beta)$}
\label{fig08}
\end{figure}

Note that the contours $\sigma_1$, $\sigma_2$, $\sigma_3$, $\sigma_4$ 
can be deformed into each other. 
Topologically, the relative positions of the contours
$\sigma_j$
and the infinity 
points $J_1, \dots , J_4$ on the torus $\gH$ are shown in Fig.~\ref{fig08}.
One can see that carrying the contours in the
positive $\alpha$ direction corresponds to moving the observation point in the $(m,n)$-plane
in the counter-clockwise (positive) direction. 
The representations are converted into each other, and every time there is a region where at least two representations are valid simultaneously. 

Representing the field in the form of an integral of a differential form 
over some compact Riemann surface puts the problem into the context of 
Abelian differentials and integrals. 
Some benefits can be gained from this. 
Namely, one can notice that, for example, (\ref{eq0107b}) is a period of 
an elliptic integral of a general form.  
One can apply the theorem by Legendre stating that a general elliptic integral 
can be represented as a linear combination of four basic elliptic integrals 
with rational coefficients \cite{Bateman1955},~p.297. Since the periods of the 
integrals are studied, the coefficients should be constant. Following the proof of the 
theorem, one can conclude that there should exist recursive relations between 
the values of $u(m,n)$ enabling one to express any $u(m,n)$ from several 
initial values computed by integration. Such a system of   
recursive relations is presented in Appendix~A. These relations can be used for 
an efficient computation of the Green's function.




\section{Diffraction by a Dirichlet half-line}


\subsection{Problem formulation}
\label{subsec:form}

\begin{figure}[ht]
\centering
\includegraphics[width=0.3\textwidth]{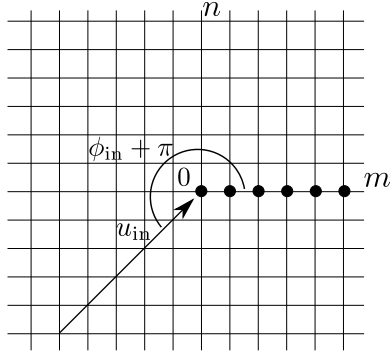}
\caption{Geometry of the problem of diffraction by a half-line. The black circles show the position of the Dirichlet half-line}
\label{fig12}
\end{figure}

Consider the following scattering problem. Let the homogeneous Helmholtz 
equation (\ref{eq0101h}) be 
satisfied by the field $u(m,n)$
for all $(m,n)$ except the half-line 
$n=0,m\geq 0$ (see Fig.~\ref{fig12}). 
On this line we impose the ``Dirichlet boundary condition''
$u = 0$. 

Let $u$ be a sum of an incident and a scattered field:
\begin{equation}
u(m,n)= u_{\rm in} (m,n) + u_{\rm sc} (m,n),
\label{uinusc}
\end{equation}
where 
\begin{equation}
\label{uin}
u_{\rm in} (m,n) = x_{\rm in}^m \, y_{\rm in}^n.
\end{equation}
The point $(x_{\rm in}, y_{\rm in})$ belongs to $\gH$, i.~e.\ 
\[
\hat D(x_{\rm in}, y_{\rm in}) = 0.
\]

We assume that the point $(x_{\rm in}, y_{\rm in})$ is taken
on the line of ``real waves'' $\beta = \pi$.  
Introduce a real parameter $\phi_{\rm in}$, which has the meaning
of angle of incidence linked with $(x_{\rm in}, y_{\rm in})$
by the relation 
\begin{equation}
\frac{y_{\rm in} - y_{\rm in}^{-1}}{x_{\rm in} - x_{\rm in}^{-1}}
= \tan \phi_{\rm in} .
\label{phiin}
\end{equation}
Angle $\phi_{\rm in}$ is the angle of propagation of the 
incident wave, while the angle of incidence is, obviously, 
$\phi_{\rm in} + \pi$. 
Let be $-\pi/ 2< \phi_{\rm in} < \pi/2$, and, thus, 
\begin{equation}
 |x_{\rm in}| < 1.
\label{eq:incond}
\end{equation}
Equation (\ref{phiin}) with condition (\ref{eq:incond})
defines two point, and only one of them belongs to the line 
$\beta = \pi$. We remind that this line is the branch of the set 
defined by (\ref{rewaves}) tending to a part of the circle $|x| = 1$
as ${\rm Im}[K] \to 0$. 
By construction of
the coordinates $(\alpha, \beta)$,
the coordinate $\alpha$ of the point $(x_{\rm in}, y_{\rm in})$
is equal to $\phi_{\rm in}$.
 

The scattered field should obey the limiting absorption principle, i.~e.\
decay as $\sqrt{m^2 + n^2} \to \infty$. 

The problem formulated here can be solved using the Wiener--Hopf method \cite{Eatwell1982, Slepyan1982, Sharma2015b}. 
This solution can be found in Appendix~B. Here, however, our aim is to develop 
the Sommerfeld integral approach for this problem.

\subsection{Formulation on a branched surface}
\label{SecBranch}

Here, following A.~Sommerfeld, we use the principle of reflections 
to get rid of the scatterer and, instead, to formulate a propagation 
problem on a coordinate plane with branching. 

Parametrize the points $(m,n)$ by the coordinates 
$(r, \phi)$:
\begin{equation}
m = r \cos \phi, \qquad n = r \sin \phi. 
\label{eq0401}
\end{equation}
Indeed $(r, \phi)$ take values from a discrete set. 

\begin{figure}[ht]
\centerline{\epsfig{file=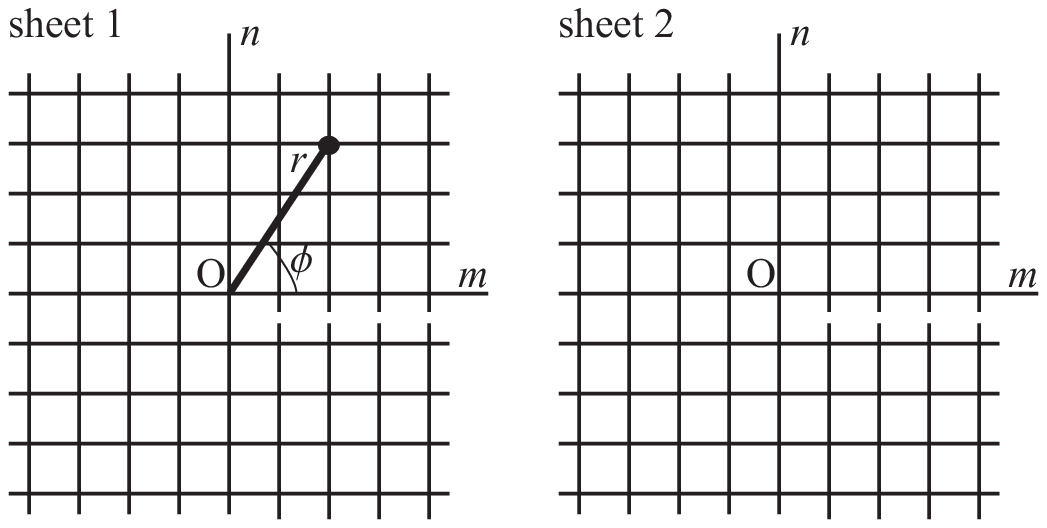,width = 9cm}}
\caption{Scheme of dicrete branched surface $\gS_2$}
\label{fig_lattice}
\end{figure}

Initially, $\phi$ belongs to $\gC$. 
Allow $\phi$ be $4\pi$-periodic,
i.~e.\ 
let $\phi$ and $\phi + 4 \pi$ mean the same, but let 
$\phi$ and $\phi + 2\pi$ correspond to different points.
Taking the points $(r, \phi)$ with such $\phi$
allows one to construct 
a branched planar lattice.
The scheme of this lattice is shown in Fig.~\ref{fig_lattice}. 
The lattice is composed of two discrete sheets. The origin $O=(0,0)$
is common for both sheets. The nodes $n = 0$, $m > 0$ of the first sheet 
are linked 
with corresponding nodes $n = -1$, $m > 0$ of the second sheet.  
The nodes $n = 0$, $m > 0$ of the second sheet 
are linked 
with corresponding nodes $n = -1$, $m > 0$ of the first sheet.  
The resulting discrete branched surface is referred to as $\gS_2$
hereafter. 

Each point of $\gS_2$ except $(0,0)$ has exactly four neighbors
in the lattice. Thus, one can look for a 
function $\tilde u$ defined on $\gS_2$ and obeying equation 
(\ref{eq0101h}) on $\gS_2 \setminus O$.

Let $u$ be the solution of the diffraction problem 
formulated in the Section~\ref{subsec:form}. 
Define the function $\tilde u$ on $\gS_2$ 
by the following formulae:
\begin{equation}
\tilde u(r, \phi) = u(r , \phi) \qquad \mbox{ for } 0 \le \phi \le 2\pi, 
\label{eq1201a}
\end{equation}
\begin{equation}
\tilde u(r, 4\pi - \phi) = -u(r , \phi) \qquad \mbox{ for } 2\pi \le \phi \le 4\pi. 
\label{eq1201b}
\end{equation}
We assume also that $\tilde u (0,0) =0$.

One can check that $\tilde u$ obeys equation (\ref{eq0101h}) on 
$\gS_2 \setminus O$. This statement is trivial for $\phi$ not equal 
to $0$ or $2\pi$ and follows from the symmetry of the field for 
$\phi = 0, 2\pi$.  

The new field $\tilde u$ has two incident field contributions. They are the 
plane wave $x_{\rm in}^m \, y_{\rm in}^n$ on sheet~1, 
and the plane wave $-x_{\rm in}^m \, y_{\rm in}^{-n}$
on sheet~2. The second wave is the reflection of the first one in the 
``mirror'' coinciding with the half-line $\phi = 0$.
Note that $(x_{\rm in}, y_{\rm in}^{-1}) \in \gH$.

We can now formulate the diffraction problem on the 
branched surface:

{\em Find a field $\tilde u$ defined on $\gS_2$, 
obeying equation (\ref{eq0101h}) on $\gS_2 \setminus O$, and equal to zero at the origin. 
The difference $\tilde u - u_{\rm in}$ with 
\[
u_{\rm in}(m,n) = w_{m,n}(x_{\rm in}, y_{\rm in})
\quad \mbox{ for } \quad 
0 \le \phi \le 2\pi,
\]
\[
u_{\rm in}(m,n) =-w_{m,n}(x_{\rm in}, y_{\rm in}^{-1})
\quad \mbox{ for } \quad
2\pi  \le \phi \le 4\pi.
\] 
should decay exponentially as $r \to \infty$.
}


\subsection{The structure of the Sommerfeld integral for the half-plane problem}
\label{Sec:3.3}

Being inspired by the classic Sommerfeld integral 
(Appendix~C), we are building an analog of the Sommerfeld integral 
on the surface~$\gH$ for finding the field~$\tilde u$. 

For this, consider an analytic manifold $\gH_2$ that is a two-sheet covering of 
$\gH$, such that the variable $\alpha$ takes values in $[0 , 4 \pi]$ (its period is doubled 
with respect to that on $\gH$), while $\beta$ still takes values in~$\gC$.
The manifold $\gH_2$ can be imagined as two copies of $\gH$, cut along the line 
$\alpha = 0$, put one above another, and connected in a single torus. 
We assume that all functions on $\gH_2$ are $4\pi$-periodic with respect to 
$\alpha$ and $2\pi$-periodic with respect to $\beta$. 

The covering $\gH_2$ has eight infinity points. 
Beside the points $J_1, \dots , J_4$ keeping the old coordinates $(\alpha, \beta)$, 
they are points 
\[
J'_1 \, : \, (9\pi/ 4, 0), 
\qquad 
J'_2 \, : \, (15\pi/ 4, 0), 
\qquad 
J'_3 \, : \, (13\pi/ 4, 0), 
\qquad 
J'_4 \, : \, (11\pi/ 4, 0).
\]

The Sommerfeld integral for the field on $\gS_2$ is an expression 
\begin{equation}
\tilde u(m,n) = \int_{\Gamma_j} w_{m,n}(p) \, A(p) \Psi,
\label{ZomPlane}
\end{equation}
where 
$\Psi$ is defined by (\ref{eq0301}),
$p = (x,y)$ is a point on $\gH_2$, 
$A(p)$ is the {\em Sommerfeld transformant\/}
of the field that is a function meromorphic on $\gH_2$, and thus, possibly, double-valued on 
$\gH$. We explain below that $A(p)$ should have two poles on the ``real waves'' line:
with $\alpha = \phi_{\rm in} + 2\pi$ and $\alpha = 4\pi - \phi_{\rm in}$, 
corresponding to the incident plane waves. 

An important part of the Sommerfeld method is the choice 
of the contour of integration in (\ref{ZomPlane}). 
By analogy with the continuous case, we need to construct 
a family of contours depending on the angle of observation $\phi$, i.~e.\ the contours 
$\Gamma = \Gamma(\phi)$. We consider a discrete family of contours, i.~e.\ there is 
a set $\Gamma_j$ with $j$ cyclically depending on $\phi$. In the continuous case 
(Appendix~C), 
the family $\Gamma (\phi)$ is formally continuous, 
but it also can be reduced to its finite subset. 

Let such a family of contours $\Gamma_j$ be constructed.
Then
for each point of $\gS_2$ and each contour $\Gamma_j$ it is possible to say, whether
the field at the point is described by the integral (\ref{ZomPlane}) with this 
contour $\Gamma_j$ or not. We will say that the point is described by the contour 
if the answer is ``yes''.

Formally, the family $\Gamma_j$ should have the  
following properties: 
\begin{enumerate}
\item
For each point $P$ of $\gS_2 \setminus O$ there should be a contour 
$\Gamma_j$ describing the point $P$ and at four its neighbors in the lattice. 

\item 
If a point $P \in \gS_2 \setminus O$
with coordinates $m,n$ is described  by two contours $\Gamma_j$ and 
$\Gamma_{j'}$ then $\Gamma_{j'}$ should be transformable 
into $\Gamma_j$ without crossing the singularities of the integrand 
of (\ref{ZomPlane}) taken for given $m,n$.  
\end{enumerate}

The second condition states that the field is describes consistently, while 
the first condition states that the field obeys 
$(\ref{eq0101h})$.

Note that in the previous section, a system of contours is built 
for the plane wave decomposition (\ref{eq0303}). 
The contours are 
$\Gamma_1^\sigma = \sigma_1$, 
$\Gamma_2^\sigma = \sigma_4$, 
$\Gamma_3^\sigma = \sigma_3$, 
$\Gamma_4^\sigma = \sigma_2$. 
they are cyclically changed as $\phi$ increases.
However, we cannot use this 
(or similar) family now, since the transformant $A(p)$ has poles on the line 
$\beta = \pi$ corresponding to the incident plane waves, and thus 
the contours $\sigma_j$ cannot be transformed one into another without hitting these poles. 

That is why, we need a new family of contours $\Gamma_j$ that 
do not cross the ``real waves'' line $\beta = \pi$. We still assume that the
contours should be obtained one from another by carrying  along the $\alpha$-axis.  

\begin{figure}[ht]
\centerline{\epsfig{file=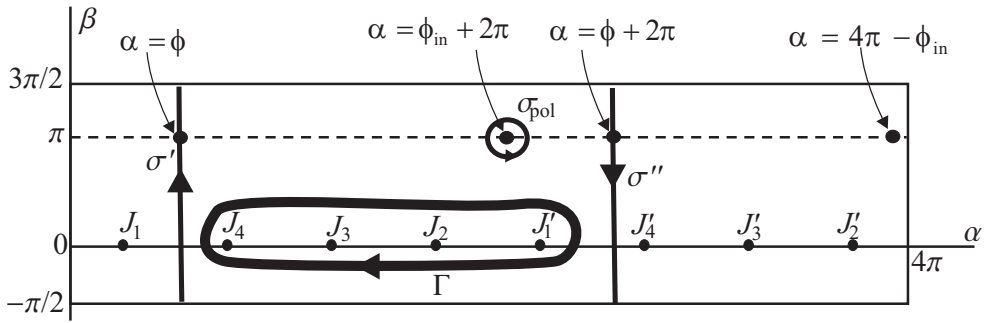,width = 13cm}}
\caption{Contour  $\Gamma$ 
and its deformation into a sum of two saddle point contours and a polar contour}
\label{fig14a}
\end{figure}

Introduce contour $\Gamma$ 
on~$\gH_2$ as it is shown in Fig~\ref{fig14a}. 
The manifold $\gH_2$ is the rectangle $0 \le \alpha \le 4\pi$,
$-\pi/2 \le \beta \le 3\pi /2$. 
Contour $\Gamma$ encircles the infinity points 
$J_4$, $J_3$, $J_2$, $J_1'$.

The family $\Gamma_j$ 
is set as follows. For $j = 1 \dots 8$ declare that the set of points $(r, \phi)$
with 
\[
(j-1)\pi /2 < \phi < (j+1)\pi /2  
\]
is described by the contour $\Gamma_j = \Gamma + (j - 1)\pi /2$. 
The origin is described by all contours.


Here we introduced the notation $\Gamma + \delta$, where $\delta$ is some angle,
that should 
be understood as follows. Contour $\Gamma$ can be considered as a set of points 
$(\alpha, \beta)$ on~$\gH_2$. Contour $\Gamma + \delta$ is $\Gamma$ shifted in the 
$\alpha$ direction by $\delta$, i.~e.\ 
changed by the transformation 
\[
(\alpha , \beta) \to (\alpha+ \delta, \beta). 
\]
The direction of $\Gamma + \delta$ is kept the same as for $\Gamma$. 

Each of the contours $\Gamma + \pi l/2$ has exactly four infinity points inside. 

The check that the family $\Gamma_j$ obeys condition~1 for the integration contour 
is trivial. To provide the validity of condition~2 one, generally, has to impose 
an additional requirement on~$A$. Let $A(p)$ have no singularities except the poles at 
$(\phi_{\rm in} + 2\pi, \pi)$
and $(4\pi-\phi_{\rm in} , \pi)$.
Let for example the point $(r, \phi)$ be described by $\Gamma_1$ and $\Gamma_2$, 
i.~e.\ 
\begin{equation}
\pi/2  \le \phi \le \pi.
\label{range}
\end{equation}
Let us show that the contour $\Gamma_1$ can be safely transformed into $\Gamma_2$. 
For this, we need to show that the integrand is non-singular at the infinity points 
$J_4$ and~$J_4'$.
Function $A$ is analytic there by our assumption. The form $\Psi$ is analytic everywhere on 
$\gH_2$. Consider the function $w_{m,n}$. Infinities $J_4$ and $J_4'$ correspond to 
$x \to  \infty$, $y \to 0$, and the range (\ref{range}) means that $m \le 0$, $n \ge 0$. 
Under these conditions, $w_{m,n} = x^m y^n$ is analytic. 

Any other pair of neighboring contours $\Gamma_j$ is checked the same way. 
Finally, the set of contours $\Gamma_j$ obey both conditions imposed on the 
Sommerfeld integral contours.


\subsection{Properties of the Sommerfeld integral}
\label{Sec:3.4}

Let be $0 \le \phi \le \pi$, so 
the contour $\Gamma_1 = \Gamma$ can be used in (\ref{ZomPlane}). 
Consider the far field, i.~e.\
build the asymptotics of $u(mN, nN)$
as $N \to \infty$ for fixed $m,n$.

As usual, such an asymptotics can be built by applying the saddle-point 
method (we follow~\cite{Martin2006}).
Let us find the saddle points on~$\gH_2$.
Represent $w_{mN, nN} (x, y)$ as 
\[
w_{mN, nN} (x, y) = 
\exp \{ 
i N (m \xi_1 + n \xi_2 )
\}
=
\exp \{ 
N (m \log(x) + n \log(\Xi(x)))
\}
.
\]
A saddle point $x_*$ corresponds to 
\[
\frac{d}{dx} \left(
m \log(x) + n \log(\Xi(x))
\right) = 0,
\]  
i.~e. 
\begin{equation}
\frac{d\Xi(x_*)}{dx} = - \frac{m}{n}\frac{\Xi(x_*)}{x_*}
\label{eq:sadpoint}
\end{equation}
This equation can be solved explicitly, since $\Xi(x)$ is given by 
(\ref{eq0111b}). 
As the result, we get a multivalued expression $x_* = x_* (m/n)$. We do not put this formula here 
due to its unwieldiness. 

Using (\ref{eq0116}), one can rewrite equation 
(\ref{eq:sadpoint}) as
\begin{equation}
\frac{y_* - y_*^{-1}}{x_* - x_*^{-1}} = \frac{n}{m}, \qquad y_* = \Xi(x_*).
\label{eq3100}
\end{equation}
Indeed,
\[
\frac{n}{m} = \tan \phi, 
\]
and the set of the points corresponding to the solutions 
of (\ref{eq3100}) for real $\phi$ is the line of ``real waves''
$\beta = \pi$, introduced above. 

Note that formally the points on the line $\beta = 0$ also satisfy the equation 
(\ref{eq3100}), but such saddle points yield asymptotic 
terms that are exponentially small. 

There are two values of $x_*(m/n)$ on $\gH$
belonging to the ``real waves'' line.
They are the points $(\pm \arctan(n/m) , \pi)$.
One of them corresponds to the wave going {\em from\/} the origin, 
and the other corresponds to the wave coming {\em to\/} the origin. 
Indeed, on $\gH_2$ there are two values $x_*(m/n)$ corresponding to the 
waves going from the origin. 
Denote these points of $\gH_2$ by $p'(m/n)$ and $p''(m/n)$. They have
$(\alpha, \beta)$-coordinates $(\arctan(n/m), \pi)$
and $(\arctan(n/m) + 2\pi, \pi)$.

Deform 
contour $\Gamma$ as it is shown in Fig.~\ref{fig14a}.
Namely, the deformed contour consists of two saddle-point loops 
$\sigma' (m/n)$ and $\sigma'' (m/n)$ passing through the aforementioned 
saddle points, and the loop $\sigma_{\rm pol}$ encircling the poles other than the 
infinities. 
Only the poles located between $\sigma' (m/n)$ and $\sigma'' (m/n)$ 
fall inside~$\sigma_{\rm pol}$. 
On the ``real waves'' line, the poles that fall between the  
contours 
$\sigma' (m/n)$ and $\sigma'' (m/n)$ should
have coordinate $\alpha$ falling in the range 
\[
\phi < \alpha < \phi+2\pi.
\]

The standard procedure of the saddle-point integration gives the cylindrical wave.
The directivity of this wave is proportional to $A(p'(m/n)) - A(p''(m/n))$, 
and this combination seems typical for a Sommerfeld integral.
The polar terms provide the incident plane waves 
(the initial wave $w_{m,n} (x_{\rm in}, y_{\rm in})$
and its mirror image) in the regions of their geometrical visibility. 
The initial plane wave corresponding to the pole 
$\phi_{\rm in} + 2\pi$ is visible in the range 
\[
\phi_{\rm in} < \phi < \phi_{\rm in} + 2\pi, 
\]
while the reflected wave corresponding to the pole 
$4\pi - \phi_{\rm in}$ is visible in the range 
\[
 2\pi - \phi_{\rm in} < \phi < 4\pi - \phi_{\rm in} .
\]

\subsection{Functional problem for the transformant $A(p)$ and its solution}
\label{Sec:3.5}

Let us formulate the functional problem for $A(p)$:

{\em
The transformant $A(p)$ should be meromorphic on $\gH_2$, that is a two-sheet covering of 
$\gH$ introduced above.
$A(p)$ can only have simple poles at two 
points of $\gH_2$, corresponding to the incident waves (the initial wave and its 
mirror reflection).
They are the points $(\phi_{\rm in} + 2\pi, \pi)$
and $(4\pi - \phi_{\rm in} , \pi)$ in the $(\alpha, \beta)$-coordinates.
The residues of the form $A \Psi$ at these points should be equal to 
$-(2\pi i)^{-1}$ and $(2\pi i)^{-1}$, respectively. 
} 

It should be clear from the consideration above, that these conditions are sufficient 
to make the Sommerfeld integral (\ref{ZomPlane}) describe a solution of the 
diffraction problem formulated above for~$\gS_2$.

\begin{figure}[ht]
\centering
\includegraphics[width=0.5\textwidth]{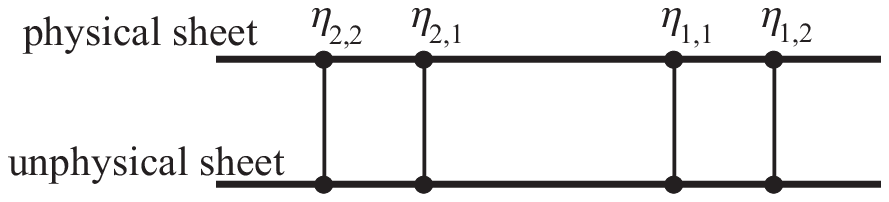}
\caption{Scheme of Riemann surface $\gR$}
\label{fig13b}
\end{figure}

Let us reformulate the problem for $A(p)$ in more usual terms. 
Consider $x$ as the main variable.
As we mentioned above, the torus $\gH$ is the preimage of the Riemann surface 
$\gR$ shown in Fig.~\ref{fig03}.  
A usual scheme \cite{Alekseev2004} of this surface is shown in Fig.~\ref{fig13b}. 
The horizontal lines are the samples of the complex planes of variable~$x$. The circles are the 
branch points. The vertical lines are the connections between sheets due to the 
branching. This scheme should be completed by the scheme of the cuts drawn in the complex 
$x$-plane shown in Fig.~\ref{fig03}. 
 
\begin{figure}[ht]
\centering
\includegraphics[width=0.5\textwidth]{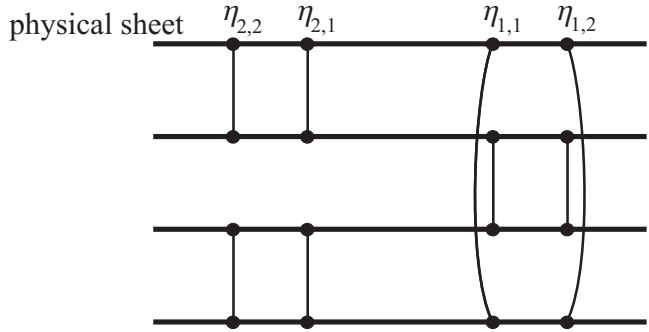}
\caption{Scheme of Riemann surface $\gR_2$}
\label{fig14b}
\end{figure}

One can see that $\gR$ is the covering over $\overline{\mathbb{C}}$ with branching. 
In the same way we can say that $\gH_2$ is the preimage of some Riemann surface 
$\gR_2$ over $\overline{\mathbb{C}}$: 
$
\gR_2 = \zeta (\gH_2)$. 
The scheme of $\gR_2$ over $\overline{\mathbb{C}}$ is shown in Fig.~\ref{fig14b}. 
The coverings can be described 
by the diagram
\[
\gR_2 \to \gR \to \overline{\mathbb{C}},
\]
where the second mapping is $\zeta$. The first mapping has no branch
points. 

Consider the points of $\gH_2$ with $x = x_{\rm in}$. 
There are 4 such points: 
\[
p_1 :\,  (\phi_{\rm in} + 2\pi  , \pi),
\qquad 
p_2 :\,  (\phi_{\rm in} , \pi),
\qquad
p_3 :\,  (2\pi- \phi_{\rm in} , \pi)
\qquad 
p_4 :\,  (4\pi - \phi_{\rm in} , \pi), 
\]
in the $(\alpha, \beta)$-coordinates. 

 Two of these points, 
$p_1$ and $p_4$ correspond to the incident plane waves mentioned 
in the condition of the functional problem. 
Indeed, point $p_1$ has $y = y_{\rm in}$, while $p_4$ has $y = y_{\rm in}^{-1}$. 
The function $A$ should have poles at these points. 
The prescribed residues are 
\begin{equation}
{\rm Res}[A \Psi , p_1] = - {\rm Res}[A \Psi , p_4] = -(2\pi i)^{-1}.
\label{system1}
\end{equation}
According to the functional problem, $A$ should be regular at $p_2$ and $p_3$, 
otherwise the field would have the plane wave terms other than the incident components listed in the problem formulation. Thus, 
\begin{equation}
{\rm Res}[A \Psi , p_2] =  {\rm Res}[A \Psi , p_3] = 0.
\label{system2}
\end{equation}


The problem for $A(p)$, $p \in \gH$ becomes reformulated as a 
problem of finding a 
4-valued function 
$A(x)$ on $\overline{\mathbb{C}}$ having given branch points 
$\eta_{j,k}$,
single valued on a prescribed Riemann surface $\gR_2$,
and having given poles $x_{\rm in}$ with prescribed residues 
on each sheet.   This is a problem of theory of functions. 
It is easy to show that $A(x)$ should be an
{\em algebraic\/} function.
To build this function, below we guess four basis functions
having all possible symmetries on $\gR_2$ and construct
$A(x)$ explicitly using these basis functions. In general, however, 
and in particular for the right angle problem studied in the next 
section of the paper, similar problems seem more complicated and require 
application of some advanced methods.   
  
An elementary symmetrization argument shows that {\em any\/} 
function meromorphic on $\gR_2$, e.g.\ $A(x)$, 
can be written as a linear combination of four 
functions having different types of symmetry with respect to 
the substitutions of sheets of $\gR_2$.
These four functions are 
\begin{equation} 
f_0 = 1,
\qquad 
f_1 (x) = \sqrt{(x-\eta_{1,1})(x-\eta_{1,2})(x-\eta_{2,1})(x-\eta_{2,2})},
\label{eq1205}
\end{equation}
\begin{equation} 
f_2 (x) = \sqrt{(x-\eta_{2,1})(x-\eta_{2,2})},
\qquad 
f_3 (x) = \sqrt{(x-\eta_{1,1})(x-\eta_{1,2})}
\label{eq1206}
\end{equation}
(note that these functions are single-valued on~$\gR_2$).
The coefficients of the linear combinations are rational functions of~$x$.

Thus,  
\begin{equation}
A(x) = 
R_0(x) + R_1(x) f_1(x) + R_2(x) f_2(x) + R_3(x) f_3(x). 
\label{eq1207}
\end{equation} 
Our aim is to find the functions $R_j (x)$.

Let the residues of $R_j (x)$ at $x = x_{\rm in}$ be equal to 
$c_j$. Using these coefficient, find the residues of the form $A \Psi$
at the points $p_j$ on $\gH_2$: 
\begin{eqnarray*}
{\rm Res}[A \Psi , p_1] &=& (c_0 + c_1 g_1 + c_2 g_2 + c_3 g_3) / g_1,
\\
{\rm Res}[A \Psi , p_2] &=& (c_0 + c_1 g_1 - c_2 g_2 - c_3 g_3) / g_1,
\\
{\rm Res}[A \Psi , p_3] &=& -(c_0 - c_1 g_1 - c_2 g_2 + c_3 g_3) / g_1,
\\
{\rm Res}[A \Psi , p_4] &=& -(c_0 - c_1 g_1 + c_2 g_2 - c_3 g_3) / g_1,
\end{eqnarray*}
where
\[
g_1 =  f_1 (p_1),
\quad
g_2 =  f_2 (p_1),
\quad 
g_3 =  f_3 (p_1). 
\]
Taking the value $f_j (p_1)$ means that the one should take the value 
$f_j (x_{\rm in})$ on the sheet of $\gR_2$ where the point $p_1$ is located.

Solving equations (\ref{system1}), (\ref{system2}), obtain
\begin{equation}
c_0 = -\frac{1}{4\pi i} f_1 (p_1), 
\qquad 
c_2 = -\frac{1}{4\pi i} f_3 (p_1), 
\qquad 
c_1 = c_3 = 0.
\label{eq1208}
\end{equation}

Finally, we can construct $A(x)$ taking the simplest rational functions 
having prescribed 
residues at $x= x_{\rm in}$:
\[
R_j = \frac{c_j}{x - x_{\rm in}}, 
\]
and
\begin{equation}
A(p) =-\frac{ f_1 (p_1) + f_3(p_1) \, f_2(p)  }{4\pi i (x - x_{\rm in})}
.
\label{eq1209}
\end{equation}
One can check that this transformant obeys all conditions imposed on $A$.

Note that the transformant $A(p)$ has no singularities at infinity. The reason of this is 
that the Dirichlet problem is studied, and the reflected plane wave has the sign opposite to
 the incident plane wave. Thus, the sum of residues of the transformant corresponding to 
 the plane waves is equal to zero and there is no need to additional poles at the infinity points.


\subsection{Analysis of the solution obtained by the Sommerfeld integral}
\label{Sec:3.6}
Let us check directly that the Sommerfeld solution  (\ref{ZomPlane}) with the transformant (\ref{eq1209}) coincides with 
the Wiener--Hopf solution (\ref{WHsol}) (the details are given in Appendix~B). For definiteness, let us study the case $0 \le \phi < \pi$, i.~e.\ 
consider the Sommerfeld integral with contour $\Gamma$. 
Deforming contour $\Gamma$ as shown in Fig.~\ref{fig14a}, obtain:
\begin{equation}
\tilde u(m,n) = u_{\rm in}(m,n) + \int_{\sigma' + \sigma''}w_{m,n}A(p)\Psi.
\end{equation}  
Taking into account the symmetry of function $A$ as $\alpha \to \alpha + 4\pi$, 
obtain for the second term
\begin{equation}
-\frac{1}{4\pi i}\int_{\sigma'}\frac{x^{m}y^{n}}{x(y - y^{-1})}\left(\frac{f_1(p_1) + f_3(p_1)f_2(x)}{x-x_{\rm in}} - \frac{f_1(p_1) - f_3(p_1)f_2(x)}{x-x_{\rm in}} \right)dx  = 
\end{equation}
$$
= -\frac{1}{2\pi i}\int_{\sigma'}\frac{x^{m}y^{n}}{x(y - y^{-1})}\frac{f_3(p_1)f_2(x)}{x-x_{\rm in}}dx.
$$
Here we used explicit expressions for $A,\Psi$ and $w_{m,n}$.  Taking into account (\ref{eq0111w}), obtain (\ref{WHsol}).

While representation (\ref{WHsol}) seems to be simpler, the Sommerfeld integral representation (\ref{ZomPlane}) may be more suitable for computations, 
since it is an integral over several polar singularities
possibly located at the infinity points $J_1', J_2, J_3, J_4$.   Thus, the integral can be calculated using the residue theorem. This means essentially that we put our problem into the 
context of the generating functions \cite{Lando2003}.

Let us for example calculate the integral on the half-line $n=0, m>0$ where the 
Dirichlet boundary conditions should be satisfied. In this case it has poles at the points $J_3$ and $J_4$. The poles are of order $m$. Straightforward calculations show that
\begin{equation}
{\rm Res}[w_{m,n}A\Psi,J_3] = - {\rm Res}[w_{m,n}A\Psi,J_4] = -\frac{1}{(m-1)!}\lim\limits_{\tau \to 0}\frac{d^{m-1}}{d\tau^{m-1}}\left[\frac{f_1(x_{\rm in})+f_3(x_{\rm in})f_2(\tau^{-1})}{4\pi i (\tau^{-1}- x_{\rm in})(y-y^{-1})\tau}\right]
\end{equation} 
and
\begin{equation}
\tilde u(m,0) = 0, \quad m>0, 
\end{equation}   
i.~e.\ the Dirichlet boundary conditions are satisfied.

Values of $\tilde u$ on the whole grid can be calculated in a similar way. If 
$m > 0, n>0, m>n$ the Sommerfeld integral has poles at $J_3, J_4$ of orders 
$m-n$ and $m+n$, respectively:
$$
\tilde u(m,n) = -\frac{1}{(m+n-1)!}\frac{d^{m+n-1}}{d\tau^{m+n-1}}\left[\frac{\tau^{n}y^n(f_1(x_{\rm in})+f_3(x_{\rm in})f_2(\tau^{-1}))}{2 (\tau^{-1}- x_{\rm in})(y-y^{-1})\tau}\right](\tau=0)
$$
\begin{equation}
\label{Resc1}
 - \frac{1}{(m-n-1)!}\frac{d^{m-n-1}}{d\tau^{m-n-1}}\left[\frac{\tau^{-n}y^{-n}(f_1(x_{\rm in})+f_3(x_{\rm in})f_2(\tau^{-1}))}{2 (\tau^{-1}- x_{\rm in})(y-y^{-1})\tau}\right](\tau=0). 
\end{equation}
If $n>0, n>m, m>-n$ the Sommerfeld integral has poles at $J_2, J_3$ 
of orders $n-m$ and $m+n$, respectively: 
$$
\tilde u(m,n) = -\frac{1}{(m+n-1)!}\frac{d^{m+n-1}}{d\tau^{m+n-1}}\left[\frac{\tau^{n}y^n(f_1(x_{\rm in})+f_3(x_{\rm in})f_2(\tau^{-1}))}{2(\tau^{-1}- x_{\rm in})(y-y^{-1})\tau}\right](\tau=0)
$$
\begin{equation}
\label{Resc2}
 - \frac{1}{(n-m-1)!}\frac{d^{n-m-1}}{dx^{n-m-1}}\left[\frac{x^{n}y^{n}(f_1(x_{\rm in})+f_3(x_{\rm in})f_2(x))}{2 (x- x_{\rm in})(y-y^{-1})x}\right](x=0). 
\end{equation}
Following this procedure one can obtain explicit expressions for the field in the remaining nodes.

Thus, using formulae (\ref{Resc1}),(\ref{Resc2}) 
one can  represent the field in each node $(m,n)$ as an explicit algebraic function of variables $(K, x_{\rm in})$.

Let us check that the homogeneous equation (\ref{eq0101h}) is satisfied at node $(-1,0)$.
For this, 
calculate values $\tilde u(m,n)$ 
explicitly at the nodes $(-2,0)$, $(-1,1)$, $(-1,-1)$, $(-1,0)$. 
We have:
\begin{equation}
\tilde u(-2,0) = \frac{f_1(x_{\rm in})(1- (k^2 - 4)x_{\rm in})}{x^2_{\rm in}},
\end{equation}
\begin{equation}
\tilde u(-1,1) = \frac{-2f_1(x_{\rm in})-2f_3(x_{\rm in})+(\eta_{2,1}+\eta_{2,2})f_3(x_{\rm in})}{4x^2_{\rm in}},
\end{equation}
\begin{equation}
\tilde u(-1,-1) = -\frac{2f_1(x_{\rm in})-2f_3(x_{\rm in})+(\eta_{2,1}+\eta_{2,2})f_3(x_{\rm in})}{4x^2_{\rm in}},
\end{equation}
\begin{equation}
\tilde u(-1,0) = \frac{f_1(x_{\rm in})}{x_{\rm in}}.
\end{equation}
Substituting the latter in (\ref{eq0101h}) and taking into account that $\tilde u(0,0) =0$, we see that the homogeneous Helmholtz equation is satisfied.


\section{Diffracton by a Dirichlet right angle}  
\subsection{Problem formulation}
\label{wedge1}
Let the homogeneous discrete Helmholtz equation (\ref{eq0101h}) 
be satisfied   everywhere except the domain $n \geq 0,m\geq 0$, 
which is the scatterer in this case. 
The Dirichlet conditions are set on the boundary of the scatterer:
\begin{equation}
\label{bound_condw}
u(m,0) = 0,\quad m\geq0, \qquad u(0,n) = 0,\quad n\geq0.
\end{equation}

\begin{figure}[ht]
\centering
\includegraphics[width=0.3\textwidth]{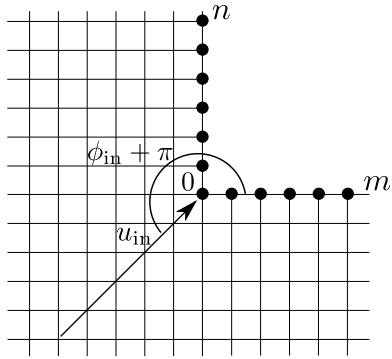}
\caption{Geometry of the problem of diffraction by an angle. Black circles show the position of the Dirichlet angle}
\label{fig12b}
\end{figure}

This corresponds to the classical 2D 
problem of diffraction by an angle (or by a wedge in 3D). The geometry of the problem is shown in Fig.~\ref{fig12b}.
The total field is a sum of the incident field  and 
the scattered field (see (\ref{uinusc})).
The incident field is the plane wave (\ref{uin}).
As before, we take the incident plane wave 
on the ``real waves line'', i.~e.\ (\ref{phiin})
is valid. The angle of propagation of the 
incident wave $\phi_{\rm in}$ obeys 
\[
0 < \phi_{\rm in} < \pi / 2,
\] 
thus 
\[
|x_{\rm in}| < 1 ,\qquad |y_{\rm in}| < 1.
\]
The scattered field  should obey the radiation condition, i.~e.\ it should decay 
at infinity.

\subsection{Formulation on a branched surface}
\label{wedge2}
Using the same consideration as in Section~\ref{SecBranch}, let us construct a 
branched lattice $\gS_3$
with three sheets.
Use the polar coordinates $(r, \phi)$ defined as (\ref{eq0401}). 
 
The lattice $\gS_3$ is the set of all integer points $(m,n)$
written in the polar coordinates $(r, \phi)$ introduced above
with 
$0 \le \phi < 6 \pi$, implying the $6\pi$-periodicity in $\phi$ of all functions. 
Obviously, $\gS_3$ is a 3-sheet covering of $\mathbb{Z} \times \mathbb{Z}$.
The origin is common for all sheets.  
 
Consider a solution $u(m,n)$ of the diffraction problem formulated above.  
Define the function $\tilde u(r, \phi)$  
on $\gS_3$ by taking
\[
\tilde u (r, \phi) =  u (r, \phi) ,
\qquad \mbox{ for }
\pi / 2 \le \phi \le  2\pi, 
\]
\[
\tilde u (r, \phi) = - u (r, 4\pi - \phi) ,
\qquad \mbox{ for }
2 \pi  \le \phi \le  7\pi / 2, 
\]
\[
\tilde u (r, \phi) =  u (r,  \phi - 3\pi) ,
\qquad \mbox{ for }
7\pi / 2   \le \phi \le  5 \pi , 
\]
\[
\tilde u (r, \phi) =  - u (r,  7\pi - \phi ) ,
\qquad \mbox{ for }
5 \pi  \le \phi \le  13 \pi /2 . 
\]
One  can check directly that field $\tilde u$ satisfies 
equation (\ref{eq0101h}) on $\gS_3 \setminus O$, i.~e.\ the boundaries 
can be discarded in the formulation of the problem. 

The new field $\tilde u$ has four incident field contributions. Define
\[
u_{\rm in}  = w_{m,n} (x_{\rm in} , y_{\rm in}) ,
\qquad 
\mbox{ for }
\pi / 2 \le \phi \le 2\pi,
\]
\[
u_{\rm in}  = - w_{m,n} (x_{\rm in} , y_{\rm in}^{-1}) ,
\qquad 
\mbox{ for }
2\pi  \le \phi \le 7\pi/2,
\]
\[
u_{\rm in} = w_{m,n} (x_{\rm in}^{-1} , y_{\rm in}^{-1}) ,
\qquad 
\mbox{ for }
7\pi/2  \le \phi \le 5\pi,
\]
\[
u_{\rm in}  = -w_{m,n} (x_{\rm in}^{-1} , y_{\rm in}) ,
\qquad 
\mbox{ for }
5\pi  \le \phi \le 13\pi/ 2.
\]
(The incident field is discontinuous, but this does not matter, 
since we are not going to construct a scattered field alone.) 
The difference $\tilde u - u_{\rm in}$ should decay as $\sqrt{m^2 + n^2 }\to \infty$.

 We can now formulate the diffraction problem on the 
branched surface:

{\em Find a field $\tilde u$ defined on $\gS_3$, 
obeying equation (\ref{eq0101h}) everywhere except the origin, 
with $u- u_{\rm in}$ exponentially decaying as $r \to \infty$.
}


\subsection{The Sommerfeld integral and the problem for the transformant}
\label{wedge3}

\begin{figure}[ht]
\centering
\includegraphics[width=0.8\textwidth]{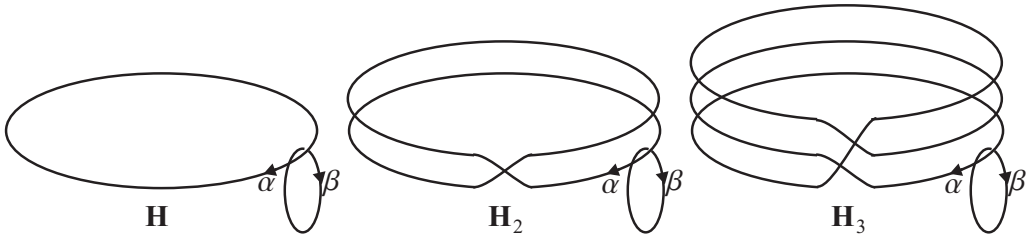}
\caption{Schemes of manifolds $\gH$, $\gH_2$, $\gH_3$ in the coordinates 
$(\alpha, \beta)$}
\label{figgHs}
\end{figure}

Analogously to the half-line diffraction problem, introduce an abstract 
analytic manifold $\gH_3$ that is a 3-sheet covering of $\gH$ without branching. 
To build this manifold, take coordinates $(\alpha, \beta)$ introduced for $\gH$
and allow $\alpha$ to take values in $[0, 6 \pi]$. This means that we take 
three copies of $\gH$ having $\alpha \in \gC$, 
all cut into ``tubes'' along the coordinate lines 
$\alpha = 0$, and attach them one to another, forming a ``ring''.
The schemes of $\gH$, $\gH_2$ and $\gH_3$ with respect to coordinates 
$(\alpha, \beta)$ are shown schematically in Fig.~\ref{figgHs}.

The covering $\gH_3$ has totally $12$ infinity points. They are the preimages 
of $x = 0$ and $x = \infty$ for the mapping $\gH_3 \to \overline{\mathbb{C}}$ 


We seek the solution $\tilde u$ in the form (\ref{ZomPlane}) with the Sommerfeld transformant $A(p)$ single-valued and meromorphic on $\gH_3$. Contours $\Gamma_j$ are 
selected in the same way as it was for half-plane problem:
They are $\Gamma_j = \Gamma + (j - 1)\pi /2$, $j = 1\dots 12$. Contour $\Gamma$
is the same as shown in Fig.~\ref{fig14a}.

The functional problem for $A(p)$ is as follows:

{\em
The transformant $A(p)$ should be meromorphic on $\gH_3$ that is a three-sheet covering of 
$\gH$.
$A(p)$ can only have simple poles at four
points of $\gH_3$, corresponding to the incident plane wave and its mirror reflections.
These poles have $(\alpha, \beta)$-coordinates
$(\phi_{\rm in} + 2\pi, \pi)$,
$(4 \pi - \phi_{\rm in}, \pi)$,
$(\phi_{\rm in} + 5\pi, \pi)$,
$(\pi - \phi_{\rm in}, \pi)$. 

The residues of the form $A \Psi$ at these points should be equal to 
$-(2\pi i)^{-1}$, $(2\pi i)^{-1}$, $-(2\pi i)^{-1}$, $(2\pi i)^{-1}$ respectively. 
}


\subsection{Solution of the functional problem for the Sommerfeld transformant}
\label{wedge4}

Unfortunately, we were not able to build a simple representation like (\ref{eq1207}) 
for $A(p)$, since we cannot guess basic algebraic functions on $\gH_3$ having 
all six types of symmetry. Here we construct $A(p)$ as an infinite series using the 
general theory of elliptic functions \cite{Gurvitz1968,Akhiezer1980,Bateman1955}.

Let $p$ denote a point on $\gH$.
Introduce a complex function $t(p)$ by the relation  
\begin{equation}
t(p) = \int^{p}_{B_1} \frac{dx}{x(y-y^{-1})} = \int^{p}_{B_1}  \Psi.
\label{elfunF}
\end{equation}
We are going to use both the function $t(p)$ and its inverse $p(t)$.
The lower limit of the integral can be arbitrary, and we take it 
equal to $B_1$ for convenience. 

Indeed, the value of $t$ depends not only on $x$ but also on the path of 
integration. 
The ambiguity is taken care about as follows. 
It is known that the mapping 
$p \to t$  maps the torus $\gH$ to a parallelogram $\gP$ in the $t$-plane with sides
$\omega_1$, $\omega_2$.     
These values are called the periods 
and are defined by the relations 
\begin{equation}
\omega_1  = \int_\sigma \frac{dx}{x(y-y^{-1})},\quad  \omega_2  = \int_\kappa \frac{dx}{x(y-y^{-1})},
\end{equation}
Contour $\sigma$ has been used above 
(this is the line $\alpha = \pi/2$ passed in the negative $\beta$-direction), 
and contour $\kappa$ is the line $\beta = 0$
passed in the positive $\alpha$-direction. Thus, function $t(p)$ is single-valued 
on $\gH$ cut along $\sigma$ and $\kappa$. 

We should make a remark that one can use the complex plane of $t$ as the coordinate plane 
for the torus instead of $(\alpha, \beta)$.  
Roughly speaking, the $\alpha$-direction is $\omega_2$, 
while the $\beta$-direction is~$-\omega_1$.


Consider a function
\begin{equation}
\label{etafunct}
E(t,\omega_1,\omega_2) = \frac{1}{t} + 
{\sum^{\infty}_{k = -\infty}}{\sum^{\infty}_{l=-\infty}}' 
\left(\frac{1}{t - \omega_1 k - \omega_2 l} + \frac{1}{\omega_1 k + \omega_2 l} + \frac{t}{(\omega_1 k + \omega_2 l)^2}\right).
\end{equation}
The term $k=0,l=0$ is excluded from the sum. This is denoted by ``prime'' notation after the sum symbol. One can show \cite{Akhiezer1980} that function $E$ is a meromorphic function of variable $t$ having poles at the points of the grid $ t  = \omega_1 k + \omega_2 l$ with  principal parts
\[
\frac{1}{t-(\omega_1 k + \omega_2 l)}.
\]

Function $E$ satisfies the following relations:
\begin{equation}
\label{etaperiod}
E(t+\omega_j,\omega_1,\omega_2) - E(t,\omega_1,\omega_2) = d_j,\quad j=1,2, 
\end{equation}
where $d_j$ are some constants.

Function $A(p)$ can be built using (\ref{etafunct}). 
Let us study $\tilde A(t) = A\left(p(t)\right)$ as a function of variable $t$. 
Note that a transition from $\gH$ to its covering $\gH_3$ looks in the $t$ coordinate
as triplicating the period $\omega_2$ (note that the $\alpha$-period is triplicated).
Namely, for $A(p)$ to be single-valued on $\gH_3$, $\tilde A(t)$ should be 
a double periodical function  with the periods $(\omega_1,3\omega_2)$.
Besides, in each elementary parallelogram the function $\tilde A(t)$
should have four poles corresponding to the incident waves.

The positions of poles on the $t$-plane are defined as follows. Let
$t_0$ be $t(x_{\rm in})$ with the integral 
in (\ref{elfunF}) taken along the shortest path along the ``real waves'' line. 
Then, one can establish the following correspondence:
\[
\phi_{\rm in} + 2\pi \longleftrightarrow  t_0 + \omega_2,  
\]
\[
4\pi - \phi_{\rm in} \longleftrightarrow   2\omega_2 - t_0,  
\]
\[
\phi_{\rm in} + 5\pi \longleftrightarrow  t_0 + 5\omega_2/2,  
\]
\[
\pi - \phi_{\rm in} \longleftrightarrow  \omega_2/2-t_0,  
\]
Thus, one can construct $\tilde A(t)$ as
\[
A\left(x(t)\right) = 2\pi i \big[-E(t-(t_0 + \omega_2),\omega_1,3\omega_2)+E(t-(2\omega_2-t_0),\omega_1,3\omega_2)-
\]
\begin{equation}
\label{Awedge}
E(t-(t_0 + 5\omega_2/2),\omega_1,3\omega_2)+E(t-(\omega_2/2-t_0),\omega_1,3\omega_2)\big].
\end{equation}
This function is triple periodic with periods 
$(\omega_1, 3\omega_2)$, since it contains two terms with $E$ and two terms 
with $-E$. As the result, the constants $d_j$ introduced
in (\ref{etaperiod}) are canceled.  
Thus, the integral (\ref{ZomPlane}) with  (\ref{Awedge}) provides the solution for the right-angled wedge problem.

Indeed, the function $A(x) = \tilde A(t(x))$ with $\tilde A$ defined by (\ref{Awedge})
is in some sense expressed in a not efficient way. It is easy to prove that 
$A(x)$ should be an algebraic function that can be expressed explicitly 
and should contain only square and cubic radicals. However, finding the 
explicit form of such an elementary function is not an elementary problem, and we prefer 
to leave the answer in elliptic functions, which are in this case most informative.  

While the obtained result seems to be cumbersome and impractical it still can be used for  calculations. Indeed,  the integral (\ref{ZomPlane}) is still a pole integral that can be calculated using residue theorem, and the formulae similar to (\ref{Resc1}) and (\ref{Resc2}) can be obtained. For this, one needs to compute the values of some elliptic functions in several specific points. The latter problem is well known \cite{Bateman1955}, and it can be solved very efficiently using the theory of $\theta$-functions.


\section{Conclusion}
\label{concl}
 
The main ideas of the paper are as follows: 

\begin{itemize}
\item
We made an elementary observation that the representation 
(\ref{eq0107b}) for the Green's function for a discrete plane 
is an elliptic integral. Thus, an application of the 
Legendre's theorem yields recursive relations, say, (\ref{eq0120}),
(\ref{eq0505}). Such recursive relations can be useful for practical 
computation of the Green's function, since they require numerical
integration only for a few first values of the field. 
We are not studying here the practical aspects of application 
of these recursive relations (basically, their stability), 
leaving this subject for another study.   

\item 
We make an observation that the plane wave decomposition 
(\ref{eq0107b}) can be considered as an integral of an analytic
differential 1-form over a contour on an analytic manifold~$\gH$.
The integral is rewritten as (\ref{eq0303}). 
The manifold $\gH$ is the set of all plane waves that can travel along the 
discrete plane, i.~e.\ it is the dispersion diagram of the system. 
Topologically, $\gH$ is a torus. 

\item
We develop the Sommerfeld integral formalism for a Dirichlet half-line 
diffraction problem. Indeed, we follow the classical continuous consideration, 
and keep the analogy as close as possible. 
Using the principle of reflections, we formulate the discrete diffraction problem
on a branched discrete surface with two sheets. Then we construct an analogue of the
Sommerfeld integral. For this, we keep the Ansatz with the differential form 
(\ref{eq0303}) (now it reads as (\ref{ZomPlane})). By analogy with the continuous case, 
we look for the Sommerfeld transformant $A(p)$ that is single-valued on a double-valued covering of $\gH$. The latter is denoted by $\gH_2$. 
The Sommerfeld contour is chosen in a way most close to the continuous case. 
We formulate and solve a functional problem for~$A(p)$. 
The latter is the classical problem of finding 
an algebraic function on a given Riemann surface having given singularities. 

\item
We make an observation that, unlike in the continuous case, the Sommerfeld integral 
has certain computational advantages {\em without\/} a transformation into the combination of the 
saddle-point contours and the plane wave residual components. Namely, the Sommerfeld integral can be taken by computing its residues at infinity points. This yields 
explicit representations
 (\ref{Resc1}) and (\ref{Resc2}), 
so numerical integration is not needed at all for the diffraction problem.
Surprisingly, the solution of the diffraction problem seems simpler than that of the 
problem for the Green's function. 

\item 
We are trying to develop the Sommerfeld formalism for another problem, to which the
reflection principle can be applied. This is the problem of diffraction by a Dirichlet 
right angle. By analogy, we write down the Sommerfeld integral and look for the
transformant $A(p)$. The transformant is an algebraic function single-valued on a 3-sheet covering of $\gH$ named $\gH_3$. Unfortunately, the problem for the transformant is more complicated in this 
case, and we have found its solution only in terms of the elliptic functions 
(\ref{Awedge}).

\end{itemize}

\bibliographystyle{unsrt}
\bibliography{Bibliography}
%
%
%
%
%
%
%




\section*{Appenidx A. Computation of the Green's function using the recursive relations}


We study here the problem for the Green's function (\ref{eq0101}).
Our aim is to find a way to tabulate function $u(m,n)$ for some set of values $m,n$. 
It is clear 
that 
\begin{equation}
\label{sym_green}
u(m,n) = u(-m,n) = u(m,-n) = u(-m,-n).
\end{equation}
Thus, one should tabulate $u(m,n)$ only for non-negative $m,n$.  

Let it be necessary to tabulate all $u(m,n)$ with
\[
|m| + |n| \le N.
\]
A naive approach requires $\sim N^2$ computations of the integral. However, 
here we show that one can compute  only two integrals, then using ``cheap'' recursive relations.
 
Perform the consideration in two steps. On the first (easy) step 
we assume all $u(m,0)$ to be known, and prove that  all  other values $u(m,n)$
can be found by recursive relations.
On the second (more complicated) step we prove that all $u(m,0)$ 
can be found from the first two values. 

\vskip 6pt
\noindent 
{\bf Step 1.}
Compute the values of $u(m,n)$ row by row. Each row is a set of values with 
$m\ge0$, $n\ge 0$, $m+ n = \mbox{const}$, i.~e.\ the rows are diagonals.

Let all values with $|m| + |n| \le M$ be already computed, and it is necessary to compute the 
values with $|m| + |n| = M + 1$. 
Then use (\ref{eq0101}) rewritten as 
a recursive relation: 
\begin{equation}
u(M+1 - n , n) + u(M+2 - n , n-1) = 
\qquad \qquad \qquad \qquad \qquad
\label{eq0120}
\end{equation}
\[
- u(M+1 - n , n-2) - u(M - n , n-1) + (4 - K^2) u(M-n, n-1).
\]
Note that all values in the right have the sum of indices $\le M$, thus they 
are computed previously. The left-hand side is a recursive relation for 
$n = 1,2,\dots, M+1$. Thus, if the values $u(m,0)$ are known, one can compute all other values by (\ref{eq0120}).

\vskip 6pt
\noindent 
{\bf Step 2.}
Let be $n = 0$. Rewrite the representation (\ref{eq0107b})
for $u(m,0)$ in the form
\begin{equation}
u(m,0) = \frac{1}{2\pi i}
\int_\sigma \frac{x^m dx}{z(x)},
\qquad 
z(x) = \sqrt{(x^2 + (K^2 - 4)x + 1)^2 - 4x^2}. 
\label{eq0501}
\end{equation}
Being inspired by the proof of Legendre's theorem for the Abelian integrals \cite{Bateman1955}, 
derive a recursive formula for $u(m,0)$. Introduce the constants $a_0, \dots, a_3$
as follows: 
\begin{equation}
a_0 = 1 , 
\qquad 
a_1 = 2 (K^2 - 4),
\qquad 
a_2 = (K^2 - 4)^2 - 2, 
\qquad 
a_3 = 2(K^2 - 4).
\label{eq0502}
\end{equation}
Using these constants one can write 
\begin{equation}
z^2(x) = x^4 + a_3 x^3 + a_2 x^2 + a_1 x + a_0. 
\label{eq0503}
\end{equation} 
Note that 
\[
\frac{d}{dx} (x^m z)  = \left( m+ 2 \right) \frac{x^{m+3}}{z(x)} 
+ 
\left( m+ 3/2 \right) a_3 \frac{x^{m+2}}{z(x)} 
+\qquad \qquad \qquad
\]
\begin{equation}
\left( m+ 1 \right) a_2 \frac{x^{m+1}}{z(x)} 
+ 
\left( m+ 1/2 \right) a_1 \frac{x^{m}}{z(x)}
+
 m a_0 \frac{x^{m-1}}{z(x)}.
\label{eq0504}
\end{equation}

Substituting this identity 
into (\ref{eq0501})
and taking into account that contour of integration 
$\sigma$ is closed, get 
\[
- \left( m+ 2 \right) u(m+3,0) 
=
\left( m+ 3/2 \right) a_3\, u(m+2,0)
+ \qquad \qquad \qquad
\]
\begin{equation}
\left( m+ 1 \right) a_2\,  u(m+1 ,0) 
+ 
\left( m+ 1/2 \right) a_1 \, u(m,0)
+
 m \,a_0 \, u(m-1,0) .
\label{eq0505}
\end{equation}
This is the recursive relation connecting five values of $u(m,0)$. Write down the latter for $m=0$  and use (\ref{sym_green}). Obtain:
\begin{equation}
\label{rec_1}
-2u(3,0) = 3/2a_3u(2,0)+a_2u(1,0)+1/2a_1u(0,0),
\end{equation}
Then, write down (\ref{eq0101}) for $m=0,n=0$ taking into account (\ref{sym_green}), and  symmetry relation $u(m,n)=u(n,m)$. Obtain:
\begin{equation}
\label{rec_3}
4u(1,0) = 1 - u(0,0)(K^2-4).
\end{equation}
Using (\ref{rec_1}-\ref{rec_3}) one can express $u(3,0)$, $u(1,0)$ in terms of $u(2,0)$ and $u(0,0)$. Then, using (\ref{eq0505}) express $u(m,0)$ in terms of $u(2,0)$ and $u(0,0)$. Thus, knowing only two values $u(2,0)$ and $u(0,0)$ it is enough to compute all other $u(m,0)$ without integration.

We should note that  similar results were obtained for the discrete Laplace equation in \cite{McCrea1940, Duffin1953,Spitzer1964,Atkinson1999,BalthVanDerPol2008}.


\section*{Appendix B. Wiener-Hopf solution for a discrete half-line problem}
Following  \cite{Sharma2015b},
let us find the solution of the half-plane diffraction problem using the Wiener-Hopf approach.  First, let us symmetrize the problem. Namely, represent the incident field (\ref{uin}) as a sum:
\begin{equation}
u_{\rm in}(m,n) = \frac{1}{2}\left(u_{\rm in}(m,n)+u_{\rm r}(m,n)\right) + \frac{1}{2}\left(u_{\rm in}(m,n)-u_{\rm r}(m,n)\right)
\equiv  u_{\rm in,s}(m,n)+u_{\rm in,a}(m,n),  
\end{equation}
where
\begin{equation}
u_{\rm in}(m,n) = x_{\rm in}^m \, y_{\rm in}^n, 
\qquad 
u_{\rm r}(m,n) = x_{\rm in}^m \, y_{\rm in}^{-n}.
\end{equation}
Then, study the equation (\ref{eq0101h}) separately for the symmetrical 
$u_{\rm sc,s}(m,n)$ 
and anti-symmetrical 
$u_{\rm sc,a}(m,n)$ part of the scattered field. Trivially, the 
anti-symmetrical scattered field  is zero:
\begin{equation}
u_{\rm sc, a}(m,n) = 0.
\end{equation}  
Thus the solution of the symmetrical problem coincide with the solution of original 
problem: 
$$
u_{\rm sc,s}(m,n) \equiv u_{\rm sc}(m,n).
$$
Without loss of generality assume that $n \geq 0$. 

Introduce direct and inverse bilateral $\mathcal{Z}$-transforms as follows:
\begin{equation}
F(z) = \mathcal{Z} \{f_n\} = \sum_{n=-\infty}^{\infty}f_nz^{-n}, 
\end{equation}
\begin{equation}
f_n =  \mathcal{Z} \{F(z)\} = \frac{1}{2\pi\i}\int_\sigma F(z)z^{n-1}dz,
\end{equation}
where $\sigma$ is a contour going along a unit circle $|z| = 1$ in the positive direction. 

Apply the $\mathcal{Z}$-transform to 
$u_{\rm sc}(m,0)$ and take into account 
the Dirichlet boundary condition
\[
u_{\rm sc}(m,0) = - u_{\rm in} (m,0) , \qquad m \ge 0.
\] 
The result is 
\begin{equation}
\label{WH1}
U_{\rm sc}(z) = -U_{\rm in}(z) + G_-(z),
\end{equation}
where 
\begin{equation}
U_{\rm in}(z) = \sum_{n=0}^{\infty}u_{\rm in} ( n,0 ) \, z^{-n} = \frac{z}{z-x_{\rm in}}
\end{equation}
and 
\begin{equation}
G_- (z) = \sum_{m=-\infty}^{-1} u_{\rm sc}(m,0) \, z^{-m},
\end{equation}
which is an unknown function analytic inside the unit circle.
Note that function $U^{\rm in}_+(z)$ is analytic outside the unit circle. 

Study a combination
\begin{equation}
w(m) = \frac{1}{2}(u_{\rm sc}(m+1,0)+u_{\rm sc}(m-1,0)) 
+ u_{\rm sc}(m, 1) + (K^2/2-2) u_{\rm sc}(m,0)
\label{defw}
\end{equation}
for  $m \in \mathbb{Z}$. 
Due to symmetry $u(m,1) = u(m,-1)$ and to equation (\ref{eq0101h}) taken for $n =0$,
$m < 0$,
\begin{equation}
w(m)=0\quad \mbox{ for } \quad m<0,
\end{equation} 
Note that the scattered field in $n>0$ consists of plane wave decaying as $n \to \infty$, 
thus
\begin{equation}
\mathcal{Z}\{ u_{\rm sc}(\cdot,n) \}(z) = 
\Xi^n(z) \,
\mathcal{Z}\{ u_{\rm sc}(\cdot,0) \}(z)  
\end{equation}
Using this relation, one can compute $\mathcal{Z}\{ u_{\rm sc}(m,1)\}$.  
Taking the $\mathcal{Z}$-transform of (\ref{defw}), obtain
\begin{equation}
\label{WH2}
W_+(z) = \frac{\sqrt{(K^2 - 4 + z + z^{-1})^2 -4}}{2}   \, U^{\rm sc}(z), 
\end{equation}
where 
\begin{equation}
W_+(z) = \sum_{m=0}^{\infty}w(m)\, z^{-m}
\end{equation} 
is analytic outside the unit circle. Combining (\ref{WH1}) and (\ref{WH2}) we obtain the following Wiener--Hopf equation:
\begin{equation}
\frac{2W_+(z)}{z \sqrt{(K^2 - 4 + z + z^{-1})^2 -4}} = -\frac{1}{z-x_{\rm in}} + \frac{G_-(z)}{z}.
\end{equation}

This equation can be easily solved. 
For this, note that the symbol can be factorized as follows: 
\[
z \sqrt{(K^2 - 4 + z + z^{-1})^2 -4} = \Upsilon(z) = 
\sqrt{(z - \eta_{2,1})(z - \eta_{2,2})}
\sqrt{(z - \eta_{1,1})(z - \eta_{1,2})} .
\] 
The first factor is analytic outside the unit circle, while the second factor is 
analytic inside the unit circle. 

Finally, the solution is as follows:
\begin{equation}
W_+(z) = - \frac{\sqrt{(z-\eta_{2,1})(z-\eta_{2,2})}\sqrt{(x_{\rm in}-\eta_{1,1})(x_{\rm in}-\eta_{1,2})}}{2(z-x_{\rm in})}.
\end{equation}
The scattered field is given by the following integral
\begin{equation}
\label{WHsol}
u^{\rm sc}_{m,n}=-\frac{1}{2\pi \i}\int_{\sigma}\frac{z^m y^{|n|}\sqrt{(x_{\rm in}-\eta_{1,1})(x_{\rm in}-\eta_{1,2})}}{\sqrt{(z-\eta_{1,1})(z-\eta_{1,2})}(z-x_{\rm in})}dz.
\end{equation}

 
\section*{Appendix C. Sommerfeld integral as an integral on the dispersion manifold} 

Let us build an analogy between the Sommerefeld integral for the continuous 
half-line diffraction problem and the Sommerfeld integral for the discrete problem.

Consider the $(x_1, x_2)$-plane, on which 
the Helmholtz equation
\begin{equation}
\Delta U(x_1,x_2) + k^2_0 U(x_1,x_2) = 0.
\end{equation}
is satisfied everywhere except half-plane
$$
x_2=0, \qquad x_1 > 0,  
$$
where Dirichlet boundary condition is imposed:
\begin{equation}
U_{\rm sc}(x_1,0) = - U_{\rm in}(x_1,0). 
\end{equation}
Here
$U_{\rm in}$ is an incident wave:
\begin{equation}
U_{\rm in} = \exp\{i k_0 (\cos \theta_{\rm in} \, x_1+
\sin \theta_{\rm in} \, x_2) \}.
\end{equation}
The radiation and Meixner conditions should be satisfied by the field in a usual way.



In the continuous case, all possible plane waves
have form 
\[
\exp \{ 
i k_0 (\cos \theta \, x_1 + \sin \theta\,  x_2)
\} 
=
\exp \{ 
i k_0 r \cos ( \phi - \theta )
\} 
\]
for the polar coordinates 
\[
x_1 = r \cos\phi, \qquad x_2 = r \sin \phi.
\] 
Thus, the set of all plane waves is the set of all (possibly complex)
angles of propagation $\theta$. This set is the strip  
\[
0 \le {\rm Re}[\theta] \le 2\pi
\]
with the edges ${\rm Re}[\theta] = 0$ and $2\pi$ attached to each other. 
Thus, topologically, the dispersion diagram is a tube.
This tube plays the role of the manifold~$\gH$ in the continuous case.
The ``real waves'' line is indeed 
the set ${\rm Im}[\theta] = 0$. 

\begin{figure}[ht]
\centering
\includegraphics[width=0.4\textwidth]{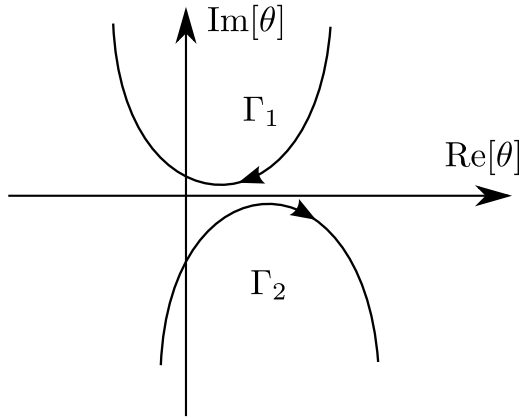}
\caption{Contours $\Gamma_1$ and $\Gamma_2$}
\label{fig11}
\end{figure}

The Sommerfeld integral has the following form:   
\begin{equation}
\label{ZomCont}
U(r, \phi) = \int_{\Gamma + \phi} A(\theta)
\exp\{ikr \cos(\theta-\phi)\} d\theta,
\end{equation}
The
contour of integration $\Gamma = \Gamma_1 + \Gamma_2$
is shown in  Fig.~\ref{fig11}.

Function $A(\theta)$ is the Sommerfeld transformant of the 
field:
\begin{equation}
A(\theta) = \frac{1}{4\pi}
\left(
\frac{\exp\{i\theta/2\}}{\exp\{i\theta/2\}-\exp\{i(\theta_{\rm in}+ \pi)/2\}}
-
\frac{\exp\{i\theta/2\}}{\exp\{i\theta/2\}-\exp\{(3\pi -\theta_{\rm in})/2\}}
\right),
\end{equation}
It is double-valued on $\gH$. Thus, one can define 
a two-sheet covering of $\gH$ named $\gH_2$, on which $A$ is single-valued.   
Such a covering is the strip 
\[
0 \le \phi \le 4 \pi 
\]
with the edges attached to each other. This covering is analogous to 
$\gH_2$ for the discrete case. 

The integral (\ref{ZomCont}) is analogous to (\ref{ZomPlane}). The exponential 
function plays the role of $w_{m,n}$, and $d\theta$ plays the role of~$\Psi$.

The transformant $A(\theta)$ has poles corresponding to the 
incident and the reflected wave. The poles belong to the 
``real line'' The contour $\Gamma + \phi$ is chosen in such a way 
that it does not hit the poles as $\phi$ changes.

\end{document}